\documentclass[letterpaper,11pt]{article}

\usepackage[margin=1in]{geometry}
\usepackage{amsmath}
\usepackage{amssymb}
\usepackage{amsfonts}
\usepackage{float}
\usepackage[font=small,labelfont=bf]{caption}
\usepackage{xfrac}
\usepackage{placeins}
\usepackage{bm}
\usepackage{amsthm}
\usepackage{xcolor}
\usepackage{diagbox}
\usepackage{arydshln}
\usepackage{subcaption}

\usepackage[ruled]{algorithm2e}
\usepackage{algorithmic}

\graphicspath{{figures/}}

\newcommand{\bb}[1]{\mathbb{#1}}
\newcommand{\mc}[1]{\mathcal{#1}}
\newcommand{\mf}[1]{\mathfrak{#1}}
\newcommand{\tc}[1]{\textcolor{#1}}
\newcommand{\trial}{\mathcal{U}}

\newcommand{\eps}{\varepsilon}
\newcommand{\p}{\partial}
\newcommand{\bs}{\boldsymbol}

\newcommand{\nE}{\hat{\bs E}}
\newcommand{\nH}{\hat{\bs H}}

\newcommand{\Nex}{\mc{N}^{\text{Yb}}_{\text{excited}}}
\newcommand{\Ngr}{\mc{N}^{\text{Yb}}_{\text{ground}}}
\newcommand{\Ntotal}{\mc{N}^{\text{Yb}}_{\text{total}}}
\newcommand{\abs}{\sigma^{\text{abs}}}
\newcommand{\ems}{\sigma^{\text{ems}}}
\newcommand{\ga}{\tilde g_a}
\newcommand{\dT}{\delta T}
\newcommand{\ndT}{\delta \hat T}
\newcommand{\bT}{\bar T}
\newcommand{\dt}{\delta t}
\newcommand{\nt}{\hat t}
\newcommand{\dn}{\delta n}
\newcommand{\bn}{\bar n}
\newcommand{\dndT}{\frac{dn}{dT}}
\newcommand{\az}{\tilde \alpha_z}
\newcommand{\Lz}{\tilde \Lambda_z}
\definecolor{dorange}{rgb}{1, 0.549, 0}

\def\div{\nabla \cdot}
\def\curl{\nabla \times}

\def\ncurl{\hat \nabla \times}
\def\hcurl{\nabla_{\hspace{-1pt}h} \times}
\def\hHcurl{H(\text{curl}, \Omega_h)}

\def\lb{\langle}
\def\rr{\rangle}

\let\oldparagraph=\paragraph
\renewcommand\paragraph[1]{\oldparagraph{#1.}}

\title{\LARGE Model and computational advancements to full vectorial Maxwell model for studying fiber amplifiers}

\author{S. Henneking$^\text{a}$, J. Grosek$^\text{b}$, L. Demkowicz$^\text{a}$\\
\small $^\text{a}$Oden Institute, The University of Texas at Austin,\\[-3pt]
\small 204 E 24th St, Austin TX 78712, USA\\
\small $^\text{b}$Directed Energy Directorate, Air Force Research Laboratory,\\[-3pt]
\small 3550 Aberdeen Ave SE, NM 87117, USA}

\date{\today}

\begin{document}

\maketitle

\begin{abstract}
We present both modeling and computational advancements to a unique three-dimensional discontinuous Petrov--Galerkin finite element model for the simulation of laser amplification in a fiber amplifier. Our model is based on the time-harmonic Maxwell equations, and it incorporates both amplification via an active dopant and thermal effects via coupling with the heat equation. As a full vectorial finite element simulation, this model distinguishes itself from other fiber amplifier models that are typically posed as an initial value problem and make significantly more approximations. Our model supports co-, counter-, and bi-directional pumping configurations, as well as inhomogeneous and anisotropic material properties. The longer-term goal of this modeling effort is to study nonlinear phenomena that prohibit achieving unprecedented power levels in fiber amplifiers, along with validating typical approximations used in lower-fidelity models. The high-fidelity simulation comes at the cost of a high-order finite element discretization with many degrees of freedom per wavelength. This is necessary to counter the effect of numerical pollution due to the high-frequency nature of the wave simulation. To make the computation more feasible, we have developed a novel longitudinal model rescaling, using artificial material parameters with the goal of preserving certain quantities of interest. Our numerical tests demonstrate the applicability and utility of this scaled model in the simulation of an ytterbium-doped, step-index fiber amplifier that experiences laser amplification and heating. We present numerical results for the nonlinear coupled Maxwell/heat model with up to 240 wavelengths.
\end{abstract}

\section{Introduction}
%
%
\paragraph{Motivation}
Optical fiber amplifiers play a critical role in our world. Optic communication technology relies on fiber amplifiers to transfer data over long distance \cite{richardson2010high}, and astronomers use fiber lasers as highly coherent light sources to calibrate telescopes \cite{wizinowich2006observatory}. Silica fiber amplifiers have emerged as excellent candidates for achieving high power outputs with great efficiency while providing superior beam quality.
The research in, fabrication of, and applications for optical fibers have greatly benefited from the growth and ubiquity of the telecommunications industry, leading to more reliable, lower cost, and higher power fiber laser systems \cite{jauregui2013fiber,grosek2017laser}.
Nonetheless, the efforts of power-scaling beam combinable fiber amplifiers have encountered numerous roadblocks in the form of nonlinear effects 
\cite{jauregui2013fiber}.
In continuous wave operation, stimulated Brillouin scattering (SBS) \cite{kobyakov2010sbs} and stimulated Raman scattering (SRS) \cite{smith1972optical} impose limits on the achievable power. One effective mitigation strategy for these nonlinear effects is to increase the core size while simultaneously decreasing the fiber length.
Unfortunately, these large-mode-area fibers permit additional guided higher-order modes, which ultimately can reduce the output beam quality if any significant amount of power is transmitted in these modes.
This happens above a certain power threshold, at which point the fundamental mode begins to exchange energy with the higher-order modes. The origins of these mode instabilities are understood to be tied to thermal effects and the interference patterns between the guided modes of the amplifier \cite{jauregui2020tmi}.


Experimental investigations of nonlinear effects in fiber amplifiers are slow, costly, and provide limited data. Indeed, the fabrication of rare-earth doped fiber amplifiers is very expensive, and manufacturing techniques limit the options and time to delivery of custom configurations. Additionally, experimental data are limited with regard to accuracy and placement of sensors. In fact, measurements are largely confined to observing the output of the fiber amplifier. Therefore, the need for modeling and simulation of fiber amplifiers is evident.

\paragraph{Methodologies}
The state of the art in numerical simulation of fiber amplifiers consists of beam propagation methods (BPMs), coupled mode theory (CMT) approaches, and a variety of other models that make certain approximations and assumptions to achieve simplified, but efficient simulations.
Some models couple to a time-dependent heat equation \cite{ward2012origin,naderi2013tmi}, others solve the thermal problem in the frequency domain \cite{hansen2012thermally,smith2013spontaneous}. Generally, these models are derived from the time-harmonic Maxwell equations and make certain assumptions to arrive at simpler models that are easier to discretize and compute. BPMs postulate that the propagating fields are guided along the longitudinal fiber direction with some propagation constant (wavenumber), and proceed to solve the corresponding field envelope by stepping along the fiber in the wave propagation direction, treating the problem as an initial value problem. Both 3D vectorial BPMs \cite{saitoh2001bpm} and, more commonly, scalar BPMs \cite{gonthier1991bpm,ward2013bpm} have been proposed. 
These BPMs tend to work better in frequency domain problems where each transverse guided mode is given its own unique frequency and wavenumber, leading to a coupled system with a different PDE for each guided mode. Otherwise, one is left with a single PDE with a given propagation constant that must capture all of the guided modes simultaneously, which strains the limits of the slowly varying envelope approximation.

Scalar models additionally assume that the propagating fields are strictly polarized in one of the transverse directions, and that the fiber is polarization maintaining, thereby eliminating two of the vector components from the equations.  Moreover, the field envelopes may be assumed to be slowly varying in the direction of propagation, reducing the model to a 2D BPM. A further simplification is made by the CMT approach that decomposes the electric field into a discrete set of propagating guided modes of the fiber, which are explicitly connected to one another via coupling coefficients \cite{naderi2013tmi, grosek2018cmt}. 

These models, posed as an initial value problem, are most amenable to forward propagating light, though they can be made to handle bi-direction light propagation; however, they tend to be very computationally intensive for such cases. For this reason, Brillouin scattering in fibers is most often modeled separately in codes written specifically for that phenomenon.

\paragraph{3D Maxwell fiber model}
This effort explains the advancements made to our unique 3D vectorial Maxwell fiber model based on a discontinuous Petrov--Galerkin (DPG) \cite{demkowicz2017dpg} finite element discretization. Rather than treating the problem as an initial value problem, this model states a boundary value problem. The goal of this approach is to make as few assumptions as possible in order to provide a tool for high-fidelity numerical simulations. The proposed model builds on the work of Nagaraj et al.\ 2019 \cite{nagaraj2018raman}, who have used a Maxwell model to simulate Raman amplification in a fiber. 
Building on that framework, we have added the ability to model the more common active gain amplification through a rare-earth, lanthanide metal dopant in the fiber core region. Like Raman gain, active gain causes our Maxwell system to become nonlinear.
Additionally, we have augmented our simulation with a thermal model that analyzes the interplay between heat deposition, the induced thermal perturbations to the refractive index of the medium, and the response of the propagating optical fields to this perturbed medium. To the best of our knowledge, this is the first fiber model that is computed using a 3D finite element discretization, and the first to simulate active gain amplification with integrated thermal response for the full vectorial Maxwell equations.

The use of the DPG method for the fiber amplifier problem is motivated by several key points \cite{henneking2020pollution}: We obtain a stable discretization that is robust with respect to the wavenumber; the ultraweak DPG Maxwell formulation exhibits numerical pollution \cite{babuska1997pollution} primarily as a diffusive effect (thereby avoiding phase errors); and, higher-order modes can be captured via local $hp$-adaptivity, driven by the DPG residual, rather than using a refined initial mesh that raises the computational cost globally.

Numerically, the proposed model poses many challenges, because it is inherently a high-frequency wave propagation problem. To counter the pollution errors effectively, the polynomial order of approximation must be raised as we increase the length of the fiber \cite{melenk2011conv,henneking2020pollution}. In fact, the number of wavelengths in a real fiber amplifier is so large that the simulation of the proposed model can only capture a tiny fraction of the entire fiber. Even much simpler models are computationally expensive, and many efforts are made to decrease the cost \cite{gopala2019equivalent}. Still, a short fiber can be used to simulate nonlinear gain, as shown in \cite{nagaraj2018raman} for a fiber of only 80 wavelengths. Drake et al.\ \cite{gopala2019equivalent} have recently investigated a scaled CMT model for active gain fiber amplifiers. The goal of our paper is to introduce the novel coupled Maxwell/heat model, and to numerically corroborate the scaling arguments with regard to nonlinear gain and heat distribution in the fiber. We emphasize that this model is unique among fiber waveguide models, and the effort described in this work strives to make progress by incorporating important physics, while maintaining computational feasibility. The investigation of relevant nonlinear phenomena (e.g., thermal lensing, SBS, SRS, Kerr effects, fiber coiling effects, etc.) in the context of the proposed model is postponed to the future.

\paragraph{Outline}
We aim to make this paper accessible to computational mathematicians as well as laser physicists and optical engineers. Thus, we will explain certain parts of our model with extra detail that may be well-known by experts in either field. In Section~\ref{sec:fiber}, we briefly introduce the general setup for an ytterbium-doped, step-index fiber amplifier. Section~\ref{sec:gain} discusses our nonlinear gain amplification model for the time-harmonic Maxwell equations, which involves two weakly-coupled systems. We make an effort to describe the non-dimensionalization in detail, because the disparate length scales in the fiber pose a modeling and computational challenge. In Section~\ref{sec:thermal}, we describe the modeling of the thermal response via coupling with the heat equation. This includes both the heat deposition in the fiber, based on the laser amplification, as well as the induced thermal polarization. The longitudinal scaling of the fiber model is discussed in Section~\ref{sec:scaling}, where we present new arguments on how to rescale the coupled Maxwell/heat model while preserving several quantities of interest. In particular, we can obtain an accurate heat distribution along the fiber amplifier for an arbitrary scaling factor. Section~\ref{sec:numerical} presents our numerical results for an ytterbium-doped active gain fiber amplifier with the full vectorial model. First, we introduce the DPG formulation of the coupled Maxwell/heat problem, and we discuss the computational complexity of the model. We then show the power distribution along the fiber, the amplifier efficiency, and convergence of the DPG residual in the nonlinear solve. We numerically corroborate the scaling arguments from the preceding section by qualitatively and quantitatively comparing the results for fibers of different length. Lastly, we show simulation results that illustrate the effect of the heating on the fiber material refractive index. In Section~\ref{sec:summary}, we summarize our work and give a brief outlook towards the wide applicability of this unique 3D model.

\section{Ytterbium-doped fiber amplifier}
\label{sec:fiber}
%
%
\paragraph{Fiber amplifiers}

\begin{figure}[ht]
	\centering
	\includegraphics[width=0.8\textwidth]{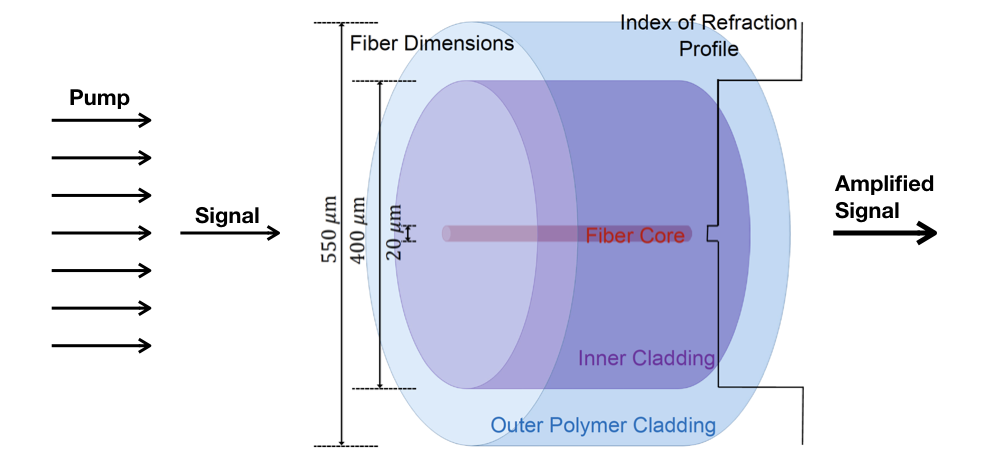}
	\caption{Schematic of a weakly-guiding, continuous-wave, double-clad, large-mode-area, step-index fiber}
	\label{fig:fiber-schematic}
\end{figure}
Figure~\ref{fig:fiber-schematic} is a schematic of a typical final-stage large-mode-area step-index fiber amplifier, not drawn to scale. The typical configuration includes a highly coherent laser signal, injected into the fiber core region, and a less coherent pump field launched into the inner cladding and core regions simultaneously (as depicted in Fig.~\ref{fig:fiber-schematic}), usually by means of the pump combiner that is spliced onto the beginning of the amplifier. However, it is also possible to core-pump the amplifier by launching the pump field into only the fiber core region as is done with the signal.  
The fact that the core region has a slightly higher refractive index than the inner cladding region, both made primarily of fused silica, allows for the signal to be guided along the fiber length in the core by total internal reflection. Thus, the fiber amplifier is a waveguide. Likewise, the polymer coating, which has a much lower refractive index than the glass, helps guide the pump light (cladding-pumped configuration) in the inner cladding region by total internal reflection. The light naturally falls into a discrete set of guided modes, where the fundamental mode has a Gaussian-like profile, which leads to the best beam quality that the fiber waveguide can output.

Active gain fiber amplifiers have doped core regions (e.g., ytterbium dopant), and the pump and signal wavelengths are chosen such that the pump field experiences very high absorption by the active dopant, while the signal field is in a regime of high emission probability of the active dopant.
Given sufficiently high pump and signal powers launched into the amplifier, a large percentage of the pump light can be converted into highly coherent laser signal light by the stimulated emission process.
This mechanism is also referred to as {\em active gain}.
Since the pump photons have a higher frequency than the signal photons, some energy is lost in this process, ultimately leading to heat deposition along the fiber via this {\em quantum defect}.
The fiber amplifier length is usually somewhere between $5\text{--}20$ m long, chosen so as to absorb a large portion of the launched pump light (e.g., $\sim\hspace*{-2pt} 95\%$), and/or to limit the onset of detrimental optical nonlinearities (e.g., SBS).
Given a strong signal seed into the fiber, spontaneous emission, amplified spontaneous emission, or even parasitic lasing are negligible issues.
Multi-photon interactions with the active dopant, e.g., photo-darkening or color-centers, or even dopant clustering, have been documented issues for some amplifiers, but can be limited by better fabrication techniques, and so will not be considered in this effort.

\section{Active gain fiber model}
\label{sec:gain}
%
%
The propagation of optical fields in fibers is governed by Maxwell's equations. We assume that the free current density and charge density vanish in the silica fiber. And, to a good approximation, optical fibers are nonmagnetic, hence we may neglect induced magnetic polarization in the material \cite{agarwal}. We assume that the signal field and the pump field each independently satisfy the time-harmonic Maxwell equations, but they are weakly coupled through nonlinear polarization terms. The time-harmonic ansatz is justified by the fact that the signal and pump fields are near monochromatic, and that the time-dependent phenomena of interest, most notably thermal effects, happen at a much slower time scale. The time-harmonic Maxwell equations in the silica fiber are,
\begin{align}
	\curl \bs E &= - i \omega \mu_0 \bs H,
	\label{eq:maxwell-3} \\
	\curl \bs H &= i \omega (\eps_0 \bs E + \bs P),
	\label{eq:maxwell-4}
\end{align}
where $\bs E$ and $\bs H$ are the electric and magnetic field vectors, respectively, $\bs P$ is the induced electric polarization, $\omega$ is the angular frequency, and $\eps_0, \mu_0$ are the electric permittivity in vacuum and the magnetic permeability in vacuum, respectively. The complex-valued time-harmonic fields $\bs E, \bs H$ are related to the real-valued time-dependent fields $\bs{\mc E}, \bs{\mc H}$ through the following ansatz:
\begin{align}
	\bs{\mc{E}}(x,y,z,t) &= \mf{Re} \left\{ \bs E(x,y,z) e^{i\omega t} \right\} , \\
	\bs{\mc{H}}(x,y,z,t) &= \mf{Re} \left\{ \bs H(x,y,z) e^{i\omega t} \right\} .
\end{align}
We assume that the fiber is aligned with its longitudinal axis centered along the $z$-axis, and the transverse coordinates are $x,y$. The core and cladding radii are denoted by $r_{\text{core}}$ and $r_{\text{clad}}$, respectively, with the corresponding material refractive indices $n_{\text{core}}$ and $n_{\text{clad}}$.

The generally nonlinear polarization term $\bs P$ may be modeled as \cite{agarwal},
\begin{equation}
	\bs P = \eps_0 \left( 
	\bs \chi^{(1)} \cdot \bs E +
	\bs \chi^{(2)} : \bs E \otimes \bs E +
	\bs \chi^{(3)}\ \vdots\ \bs E \otimes \bs E \otimes \bs E + 
	\ldots \right) ,
\end{equation}
where $\bs \chi^{(k)}$ is the $k$th-order electric susceptibility, generally given by a $(k+1)$-rank tensor. A medium with an inversion symmetry at the molecular level (e.g., $\text{SiO}_2$) exhibits only negligible second-order susceptibilities \cite{agarwal}. Therefore, in the optical fiber, $\bs \chi^{(2)} \approx \bs 0$.

In the context of fiber optics, we are interested in modeling certain nonlinear physical effects; thus it is common to write,
\begin{equation}
	\bs P = \bs P_{\text{background}} + \bs P_{\text{active gain}} + \bs P_{\text{thermal}} + \bs P_{\text{opt. nonlin.}} + \ldots
\end{equation}

The background polarization is a first-order susceptibility term and describes the real part of the refractive index. Active gain in the fiber may be modeled as a first-order susceptibility term as well, in terms of a complex perturbation to the refractive index. The thermal polarization perturbs the real part of the refractive index and may thus, too, be modeled as a first-order susceptibility. Optical nonlinearities such as SBS or SRS are third-order susceptibility terms.

\paragraph{Active gain model}

In the proposed model, we consider two electric polarization terms: the linear background polarization $\bs P_{\text{background}}$ and the polarization due to active gain $\bs P_{\text{active gain}}$. Thermal polarization will be discussed in the next section. Since there is no term inducing anisotropy in the material refractive index tensor $\bs n$ at this point, we may write the background polarization as,
\begin{equation}
	\bs P_k^{\text{background}} = \eps_0 (n^2 - 1) \bs E_k, \quad k \in \{s,p\} ,
\end{equation}
where index $k$ denotes the appropriate field, $s$ for signal or $p$ for pump. The active gain is modeled as a first-order term or complex perturbation to the refractive index \cite{nagaraj2018raman,grosek2018cmt},
\begin{equation}
	\bs P_k^{\text{active gain}} = i \eps_0 \frac{nc}{\omega_k} g_k(\bs E_{\{s,p\}}) \bs E_k ,
	\label{eq:active-gain}
\end{equation}
where the gain $g_k$ is a function of both electric fields $\bs E_{s}$ and $\bs E_{p}$, and $c$ is the speed of light in vacuum. The gain function carries units of m$^{-1}$, sometimes expressed in dB/m. A positive gain amplifies the field while a negative one causes decay. Considering these first-order susceptibilities only, the time-harmonic Maxwell equations in the fiber are:
\begin{align}
	\curl \bs E_k &= - i \omega_k \mu_0 \bs H_k ,
	\label{eq:maxwell-gain-3} \\
	\curl \bs H_k &= i \omega_k (\eps_0 n^2 \bs E_k + \bs P_k^{\text{active gain}}) .
	\label{eq:maxwell-gain-4}
\end{align}
With the gain polarization term (\ref{eq:active-gain}), equation (\ref{eq:maxwell-gain-4}) becomes,
\begin{equation}
	\curl \bs H_k = i \omega_k \eps_0 n^2 \bs E_k - \eps_0 nc g_k(\bs E_{\{s,p\}}) \bs E_k .
\end{equation}

In this effort, the gain function is expressed in terms of rate equations for a simplified version of the ion population dynamics for an ytterbium-doped fiber, as described in \cite{pask1995ytterbium}.  
It is assumed that the ytterbium (Yb) dopant is uniformly distributed throughout the fiber core.  The electrons of Yb atoms absorb and emit photons at a mean rate that is determined by experimentally measured absorption and emission cross-sections (see Table~\ref{tab:active-gain-parameters}), denoted as $\abs_k, \ems_k$ where $k \in \{ s,p \}$ (units of m$^2$/ion).  
This simplified two-manifold model only considers a single excited state and ground state for the electrons in the outer most shell of the Yb atom.  
The amplifier works by primarily absorbing pump photons, sending the electron into the excited state, and then having the majority of those excited electrons to be stimulated to emit a photon at the signal wavelength, coherent with the signal field, allowing the electron to return to its ground state.
This leads to a frequency dependent gain function as follows,
\begin{equation}
	g_k(\bs r, z, t) = 
	\ems_k \Nex (\bs r, z, t) - 
	\abs_k  \Ngr (\bs r, z, t), \quad |\bs r| < r_{\text{core}} ,
\end{equation}
where $\bs r$ represents the coordinate directions, $(x,y)$ -- Cartesian, or $(r,\theta)$ -- polar, such that $r \in [0, r_{\text{clad}}]$ and $\theta \in [0, 2 \pi)$;
and $\Ngr, \Nex$ are the ground-state and excited-state population concentrations of Yb ions expressed in ion/m$^3$. The total ion population concentration is assumed to be known, and must remain constant such that,
\begin{equation}
	\Ntotal = \Ngr + \Nex .
\end{equation}
Therefore,
\begin{equation}
	\frac{\p \Nex}{\p t} = -\frac{\p \Ngr}{\p t} .
\end{equation}
The transient equation for the excited ion population is given by,
\begin{equation}
	\frac{\p \Nex }{\p t} = 
	\sum_{k=\{s,p\}} 
	\frac{I_k}{\hbar \omega_k} (\abs_k \Ngr - \ems_k \Nex) - 
	\frac{\Ntotal}{\tau} ,
\end{equation}
where $I_k$ is the irradiance, $\tau$ is the measured Yb ion upper level radiative lifetime, $\hbar$ is the reduced Planck constant, and $I_k/(\hbar \omega_k)$ represents the photon flux.

This model, as has been done in other fiber amplifier models \cite{naderi2013tmi, ward2013bpm}, will neglect the time dynamics of the transit time of the light along the length of the fiber ($\sim\hspace*{-2pt}10$ ns), and the population dynamics of the active gain process, only considering the temporal evolution of the heat deposition and dissipation in the fiber.  The characteristic time to steady-state gain is on the order of $10 \ \mu$s, whereas the characteristic heat diffusion time is on the order of $1$ ms, both depending on the fiber configuration\footnote{Generally, higher pump and signal irradiances correspond to faster characteristic times to steady-state gain.}.  Thus, this effort will use the steady-state solution of the ion population model: 
\begin{equation}
	\bar{\mc N}_{\text{excited}}^{\text{Yb}} = 
	\frac{\sum_{k=\{s,p\}} \frac{I_k}{\hbar \omega_k} \abs_k}
	{\frac{1}{\tau} + \sum_{k=\{s,p\}} \frac{I_k}{\hbar \omega_k} (\abs_k+\ems_k)} \Ntotal ,
\end{equation}
which leads to a concise expression for the gain function,
\begin{equation}
\boxed{
	g_k \approx -\abs_k \Ntotal + (\abs_k + \ems_k) \bar{\mc N}_{\text{excited}}^{\text{Yb}}
} .
\end{equation}

\begin{table}[H]
	\caption{Active gain model: parameters}
	\centering
	\begin{tabular}{lll}
		Parameter & Value \\
		\hline
		$\abs_s$ & $6 \cdot 10^{-27}$ m$^2$/ion \\
		$\ems_s$ & $3.58 \cdot 10^{-25}$ m$^2$/ion \\
		$\abs_p$ & $1.429 \cdot 10^{-24}$ m$^2$/ion \\
		$\ems_p$ & $1.776 \cdot 10^{-24}$ m$^2$/ion \\
		$\Ntotal$ & $6 \cdot 10^{25}$ ion/m$^3$ \\
		$\tau$ & $8 \cdot 10^{-4}$ s
	\end{tabular}
	\label{tab:active-gain-parameters}
\end{table}

\paragraph{Non-dimensional Maxwell equations}
Non-dimensionalization is essential in the numerical computation of the solution to the Maxwell equations. In particular, when the scales involved are very disparate as in the case of optical fiber amplifiers. We define the non-dimensional variables $\hat x, \hat y, \hat z, \hat \omega_k, \nE_k, \nH_k$ by,
\begin{align}
	x &= l_0 \hat x,\ y = l_0 \hat y,\ z = l_0 \hat z, \\
	\omega_k &= \omega_0 \hat \omega_k , \\
	\bs E_k &= E_0 \nE_k , \\
	\bs H_k &= H_0 \nH_k ,
\end{align}
where $l_0, \omega_0, E_0,$ and $H_0$ are appropriate dimensional scales. See Table~\ref{tab:scales} for an overview of the selected dimensional scales in the fiber amplifier model. Additionally, let $\hat g_k$ denote the non-dimensional gain function with scale $g_0$.

Let,
\begin{align}
	\frac{H_0 l_0 \omega_0 \mu_0}{E_0} = 1, \quad
	\frac{E_0 l_0 \omega_0 \eps_0}{H_0} = 1;
\end{align}
then,
\[
	\omega_0 = c/l_0, \quad H_0 = \eps_0 c E_0 .
\]
The non-dimensional Maxwell system is,
\begin{align}
	\ncurl \nE_k &= - i \hat \omega_k \nH_k 
	\label{eq:maxwell-nondim-gain-3} , \\
	\ncurl \nH_k &= i n^2 \hat \omega_k \nE_k  - n l_0 g_0 \hat g_k(\nE_{\{s,p\}}) \nE_k ,
	\label{eq:maxwell-nondim-gain-4}
\end{align}
where $l_0 g_0$ is a non-dimensional quantity. The non-dimensional Amp\`{e}re--Maxwell equation (\ref{eq:maxwell-nondim-gain-4}) illustrates that for a positive gain function, the gain term can be interpreted as \emph{negative conductivity}, causing amplification/gain of the propagating field. Conversely, negative gain can be seen as positive conductivity or loss.
\begin{table}[ht]
	\caption{Dimensional scales}
	\centering
	\begin{tabular}{lll}
		Symbol & Description & Value \\
		\hline
		$l_0$ & Length & $10^{-5}$ m \\
		$\omega_0$ & Angular frequency & $\omega_0 = c/l_0$ [rad/s] \\
		$I_0$ & Irradiance & $10^{10}$ W/m$^2$ \\
		$E_0$ & Electric field & $I_0 = E_0 H_0$ [V/m] \\
		$H_0$ & Magnetic field & $H_0 = \eps_0 c E_0$ [A/m] \\
		$P_0$ & Power & $P_0 = I_0 l_0^2$ [W] \\
		$\sigma_0$ & Absorption/Emission cross-section & $10^{-26}$ m$^2$/ion \\
		$\nu_0$ & Ion population concentration & $10^{25}$ ion/m$^3$ \\
		$g_0$ & Gain & $g_0 = \sigma_0 \nu_0$ [1/m] \\
		$T_0$ & Temperature change & 1 K \\
		$t_0$ & Time & $10^{-3}$ s
	\end{tabular}
	\label{tab:scales}
\end{table}

\section{Thermal coupling}
\label{sec:thermal}
%
%
\paragraph{Thermal polarization}
First, we want to address how the heating affects the solution to the Maxwell system. We model the effect of the heating as an isotropic temperature dependence of the refractive index.
Let the ambient temperature be denoted by $\bT \equiv T_{\text{ambient}}$, and the refractive index at ambient temperature as $\bn \equiv n(\bT)$. Then, let the temperature $T$ and refractive index $n$ at any point in the fiber be given by,
\begin{align}
	T(\bs r, z, t) &= \bT + \dT(\bs r, z, t) , \\
	n(\bs r, z, t) &= \bn + \dn(\bs r, z, t) ,
\end{align}
where $\dT$ is the change in temperature, and $\dn$ the thermally induced perturbation to the material refractive index. Next, we linearize the refractive index perturbation,
\begin{align}
	n(T) &= \bn + \dndT (\bT) \dT + \frac{d^2n}{dT^2} (\bT) \frac{\dT^2}{2} + \ldots , \\
	\dn &\approx \dndT (\bT) \dT ,
\end{align}
where $dn/dT$ is the thermo-optic coefficient for silica glass, an experimentally measured value. In the fiber, we expect temperature changes of up to ca.\ $100$ K from ambient (room) temperature, a temperature regime in which the material refractive index change can be modeled with reasonable accuracy as a linear response to the temperature change. The thermo-optic coefficient is on the order of $10^{-5}\text{ K}^{-1}$ for $\text{SiO}_2$, hence one can expect induced refractive index perturbations of about three orders of magnitude smaller than the refractive index $n$ of the medium: $\dn \sim \mc{O}(10^{-3})$.

In the active gain model, the perturbed refractive index will affect the Maxwell solution through a change in the background and gain polarization:
\begin{align}
	\bs P_k^{\text{background}} &= \eps_0 (n(T)^2 - 1) \bs E_k, \\
	\bs P_k^{\text{active gain}} &= i \eps_0 \frac{n(T)c}{\omega_k} g_k \bs E_k.
\end{align}
Therefore, we can express thermal polarization explicitly as,
\begin{equation}
\begin{split}
	\bs P_k^{\text{thermal}}(T) 
	&=
	(\bs P_k^{\text{background}} (T) - \bs P_k^{\text{background}} (\bT)) \\
	&+ 
	(\bs P_k^{\text{active gain}} (T) - \bs P_k^{\text{active gain}} (\bT)) .
\end{split}
\end{equation}

\paragraph{Heat coupling model}
The thermal response in the fiber amplifier is modeled by the heat equation,
\begin{equation}
	\rho_0 C_p \frac{\p T}{\p t} - \div (\bs \kappa \nabla T) = Q ,
\end{equation}
where $\rho_0, C_p$, and $\bs \kappa$ are the mean density, specific heat, and thermal conductivity of silica glass, respectively. Appropriate boundary and initial conditions will be discussed later. We assume that the material is isotropic and homogeneous, so its thermal conductivity is uniform, i.e., $\bs \kappa(\bs r, z) = \kappa \bb I$, and all of the thermal parameters are assumed to be temperature independent. The right-hand side has the source term $Q = Q(\bs r, z,t)$ that couples the electromagnetic fields to the heat deposition in the fiber.
\paragraph{Heat deposition}
The heat source of a stimulated emission dominated amplifier can be modeled by \cite{smith2016mode,grosek2018cmt},
\begin{equation}
	Q(I_{\{s,p\}}) = - \left( g_p(I_{\{s,p\}}) I_p + g_s(I_{\{s,p\}}) I_s \right) ,
	\label{eq:thermal-loading}
\end{equation}
where $g_s$ and $g_p$ are the gain functions for signal and pump fields, respectively, and $I_{\{s,p\}}$ denotes the respective field irradiances. Therefore, $Q$ is explicitly dependent on the solution to the Maxwell equations, and thus is implicitly dependent on the temperature itself. Because gain occurs only inside the fiber core, the heat deposition will, too, only occur in the core.
\paragraph{Non-dimensional heat equation}
As discussed previously, we are interested in computing the temperature difference $\dT = T - \bT$ to obtain the perturbation to the refractive index $\dn \approx (dn/dT) \dT$. The heat equation for the temperature difference is given by,
\begin{equation}
	\rho_0 C_p \frac{\p (\dT)}{\p t} - \kappa \Delta (\dT) = Q(I_{\{s,p\}}) ,
	\label{eq:heat-temp-diff}
\end{equation}
with appropriate boundary and initial conditions.
The non-dimensional heat equation is,
\begin{equation}
	\frac{\p (\ndT)}{\p \nt} - 
	\alpha_0 \hat \Delta (\ndT) 
	= Q_0 \hat Q(\hat I_{\{s,p\}}) ,
	\label{eq:heat-nondim}
\end{equation}
where $\alpha_0$ denotes a non-dimensional diffusivity scale, and $Q_0$ a non-dimensional heat deposition scale.
\paragraph{Boundary and initial conditions}
We assume that the initial temperature in the fiber is the ambient temperature. In other words, the temperature difference $\dT$ is initially zero. At the radial inner cladding boundary with the polymer jacket we impose zero Dirichlet boundary conditions (ambient temperature), implying efficient cooling at the glass-polymer interface. In future work, it should be relatively straightforward to account for the polymer coating and more realistic heat dissipation into the ambient air and/or into the metal spool that fibers are usually coiled around.  
However, for the primary reason of keeping our computational domain smaller, this effort neglects the more realistic scenario.
At the fiber ends, we also impose homogeneous Dirichlet boundary conditions for simplicity because of the large aspect ratio between the length of the fiber and its radial width, which strongly suggests that most of the heat dissipation occurs through the radial direction rather than in the longitudinal direction.
To summarize,
\begin{align}
	\text{Initial condition: } \ndT(\bs r, z, 0) &= 0; \\
	\text{Boundary conditions: } \ndT(\bs r, z, t) &= 0 \quad  \text{ if }\
	\left\{ \begin{array}{cl}
		|\bs r| &=\ r_{\text{clad}}, \\
		z &=\ 0, \\
		z &=\ L.
	\end{array} \right.
\end{align}
\paragraph{Time-stepping scheme}
We use implicit Euler time stepping to advance the heat solution with a (dimensionless) time step $\delta \hat t$. The total time is denoted by $\hat t_{\text{max}}$. For $n$ uniform time intervals, $\delta \hat t = \hat t_{\text{max}}/n$. The implicit Euler scheme yields,
\begin{equation}
	\ndT_{n+1} - \delta \hat t \alpha_0 \hat \Delta ( \ndT_{n+1}) = 
	\ndT_n + \delta \hat t Q_0 \hat Q(\hat I_{\{s,p\}, n}) ,
\end{equation}
where the argument for the heat deposition (Maxwell fields) is taken from the previous time step; it would not be feasible to implicitly compute the nonlinear source term which itself requires solving a nonlinear Maxwell system. Because of that, the scheme is not fully implicit hence unconditional stability is not given. Instead, a suitable small time step providing stability and sufficient accuracy in time is needed. The scheme can be seen as a fully implicit Euler method with one Picard iteration. Similar time-stepping schemes have been studied in the context of DPG methods for linear problems in \cite{tofu2017time,roberts2020time}.

\begin{table}[H]
	\caption{Heat coupling model: parameters}
	\centering
	\begin{tabular}{lll}
		Parameter & Value \\
		\hline
		$C_p$ & $703$ (W$\cdot$s)/(kg$\cdot$K) \\
		$\rho_0$ & $2201$ kg/m$^3$ \\
		$\kappa$ & $1.38$ W/(m$\cdot$K) \\
		$dn/dT$ & $1.285 \cdot 10^{-5}$ 1/K
	\end{tabular}
	\label{tab:heat-coupling-parameters}
\end{table}

\section{Short fiber scaling}
\label{sec:scaling}
%
%

\paragraph{The equivalent short fiber}
A real fiber amplifier is about 5--20 meters long. The number of wavelengths in longitudinal direction is in the order of millions or tens of millions. A 3D vectorial finite element discretization with high order of approximation is computationally expensive: many degrees of freedom are required per wavelength to resolve the fields accurately. A state-of-the-art compute node is currently capable of solving the proposed model for $\mathcal{O}(10^2)$ wavelengths. Due to numerical pollution \cite{babuska1997pollution,henneking2020pollution}, the simulation of a full-length fiber is not feasible at this moment, even with a scalable parallel code. Consequently, the proposed high-fidelity model may either be viewed in the context of a multi-fidelity approach together with simplified models, or we must find a scaling argument that enables us to compute on a short fiber that preserves the physical quantities of interest from a full-length fiber with sufficient accuracy. The latter approach has been proposed for a scaled CMT model \cite{gopala2019equivalent}, since even with much simpler models, the simulation of fiber amplifiers remains computationally challenging today. Scaling only in the longitudinal direction makes physical sense for the fiber amplifier problem since the fiber waveguide performance is primarily derived from the transverse distribution of the index of refraction, and the refractive index remains relatively uniform along the fiber length. We corroborate a short fiber scaling argument through numerical results for gain polarization, and we are able to show a more rigorous mathematical argument for the scaling of thermal coupling in the fiber.

\paragraph{Gain scaling}
First, we aim to make an argument that gain can be simulated in a short fiber of a length $\tilde L$ much smaller than the real fiber length $L$. It seems natural to introduce an artificial gain scaling term proportional to $L/\tilde L$ to amplify the gain proportionally to the shortening in fiber length. Indeed, for a CMT model, it is possible to show a gain scaling argument in that way \cite{gopala2019equivalent}. In the nonlinear Maxwell problem, a convincing mathematical scaling argument may be hard to show, if possible at all, and at this point we restrict ourselves to numerical experiments. We introduce a non-dimensional short fiber gain amplifier $\ga$ in the Maxwell system:

\begin{equation} \boxed{
\begin{aligned}
	\ncurl \nE_k &= - i \hat \omega_k \nH_k \\[3pt]
	\ncurl \nH_k &= i n^2 \hat \omega_k \nE_k - n \ell_0 g_0 \ga \hat g_k(\nE_{\{s,p\}}) \nE_k
\end{aligned}
}
\label{eq:active-gain-Maxwell-model}
\end{equation}

Purely in the context of nonlinear gain, this scaling term could be viewed as an increase of the dopant population concentration $\Ntotal$ proportional to $\ga$.

\paragraph{Heat scaling}
The goal of the gain scaling is to obtain a field intensity (power distribution) in the short fiber that simulates the laser gain in a real-length fiber. The heat deposition depends only on the field intensity (cf.\ (\ref{eq:thermal-loading})). To simulate the heating in a fiber of length $L$, given a Maxwell solution on the short fiber of length $\tilde L$, a change of coordinates can be used to pull back the heat problem of the real-length fiber to a short domain. Suppose $L = 10 \text{ m} = 10^6 l_0$, and $\tilde L = 100 \ \mu\text{m} = 10 l_0$. Let,
\begin{equation}
	\hat z = \frac{L}{\tilde L} \tilde z ,
\end{equation}
so that,
\[
	0 < \hat z < 10^6\ \Leftrightarrow\ 0 < \tilde z < 10 .
\]

Denote $\az := \tilde L / L$; here, $\az =10^{-5}$. Then,
\begin{align}
	\frac{\p}{\p \hat z} &=
	\frac{\p}{\p \tilde z} \frac{\p \tilde z}{\p \hat z} =
	\az \frac{\p}{\p \tilde z}\ , \\[2pt]
	\hat \Delta &=
	\frac{\p^2}{\p \hat x^2} + \frac{\p^2}{\p \hat y^2} + \frac{\p^2}{\p \hat z^2} =
	\frac{\p^2}{\p \hat x^2} + \frac{\p^2}{\p \hat y^2} + \az^2 \frac{\p^2}{\p \tilde z^2}\ .
\end{align}
That is, the scaling yields an anisotropic diffusion operator. Equivalently, we may write,
\begin{equation}
	\hat \Delta = \tilde \nabla \cdot \left( \Lz \tilde \nabla \right) ,
\end{equation}
where
\begin{equation}
	\Lz = \left[ \begin{array}{ccc}
		1 & 0 & 0 \\
		0 & 1 & 0 \\
		0 & 0 & \az^2
	\end{array} \right]
	\text{ and }
	\tilde \nabla = \left[ \begin{array}{c}
		\sfrac{\p}{\p \hat x} \\
		\sfrac{\p}{\p \hat y} \\
		\sfrac{\p}{\p \tilde z}
	\end{array} \right] .
\end{equation}
We obtain the non-dimensional \emph{short fiber heat equation,}
\begin{equation}
\boxed{
	\frac{\p (\ndT)}{\p \nt} - 
	\alpha_0 \tilde \nabla \cdot \left( \Lz \tilde \nabla (\ndT) \right)
	= Q_0 \hat Q(\hat I_{\{s,p\}})
} .
\end{equation}

Intuitively, this makes sense because we are ``compressing'' discrete heat solution points from a real-length fiber together by a factor of $\az$ in $z$. Consequently, these points should experience very little diffusion in $z$. Another way to view this rescaling is to say that each element in $z$ is now solving the heat equation for a much longer distance in $z$, which is justified by the fact that the solution to the heat equation is very smooth. Oscillations in the temperature along $z$ may occur due to wave propagation phenomena such as transverse mode beating. But such phenomena are expected to occur at a scale linked to a certain number of wavelengths in the Maxwell solution (and thus smooth enough within each element). Therefore, the proposed short fiber heat equation is able to capture any physical heating phenomena in the real fiber due to active gain as long as enough wavelengths are computed in the Maxwell problem to exhibit the relevant wave phenomena. While the scaling argument in the Maxwell system was artificial and cannot be guaranteed to reproduce the correct physics without further investigation, the scaled heat equation does reproduce the physical results from a real-length fiber.

\section{Numerical scheme and results}
\label{sec:numerical}
%
%


\paragraph{DPG ultraweak Maxwell}
For a detailed exposition of the DPG method with optimal test functions, we refer to \cite{demkowicz2017dpg}. In our model, we use the DPG method for a conforming finite element discretization of the ultraweak Maxwell systems for the signal and pump fields. This particular weak form of the first-order system has proven superior in the context of high-frequency wave propagation problems \cite{demkowicz2012wavenumber,petrides2017multigrid}. In this section, we will omit the ``hat'' symbol for non-dimensional quantities; instead, every symbol is now understood to be non-dimensional and we will overload the ``hat'' symbol to indicate trace unknowns on the mesh skeleton. Let $\Omega$ denote the bounded fiber domain given by $\Omega := \Omega_t \times (0,L) \subset \bb{R}^3$, where $\Omega_t := \{ (x,y) : x^2 + y^2 < r_{\text{clad}}^2 \}$ is the (cross-sectional) transverse domain, and $L$ is the length of the (fiber) domain; the boundary is denoted by $\Gamma \equiv \p \Omega$; and $\Omega_h$ is a suitable finite element mesh with mesh skeleton $\Gamma_h$. The operator form of the Maxwell problem is,
\begin{equation}
\arraycolsep=2pt
\left\{
\begin{array}{rll}
	\curl \bs E_k + i \omega_k \bs H_k &= 0 &\quad \text{in } \Omega, \\
	\curl \bs H_k - i n^2 \omega_k \bs E_k + n\ell_0g_0\ga g_k \bs E_k &= 0 &\quad \text{in } \Omega, \\
	\bs n \times \bs E_k &= \bs n \times \bs E_{0,k} &\quad \text{on } \Gamma,
\end{array}
\right.
\end{equation}
where $k \in \{s,p\}$, $\bs n$ is the outward unit normal, and with appropriate boundary data. For the fiber amplifier model, we prescribe idealized perfect electrical conductor boundary conditions (i.e., vanishing tangential electric field) at the radial cladding boundary which is justified by the exponential decay of the core-guided modes in the cladding; at the input, the respective fundamental mode is prescribed to excite the waveguide; and at the fiber output, we use a stretched coordinate perfectly matched layer to avoid reflection of the propagating wave \cite{astaneh2018pml,nagaraj2018raman}. The corresponding DPG broken ultraweak formulation \cite{carstensen2016breaking} is,
\begin{equation}
\arraycolsep=2pt
\left\{
\begin{array}{rll}
	\multicolumn{2}{l}{\bs E_k, \bs H_k \in (L^2(\Omega))^3,
	\nE_k \in \hat \trial_1,
	\nH_k \in \hat \trial_2,} \\[3pt]
	(\bs E_k, \hcurl \bs F) + 
	\lb \bs n \times \nE_k, \bs F \rr_{\Gamma_h} +
	i \omega_k (\bs H_k, \bs F) &= 0,
	&\ \bs F \in \hHcurl, \\[3pt]
	(\bs H_k, \hcurl \bs G) + 
	\lb \bs n \times \nH_k, \bs G \rr_{\Gamma_h}
	- i \omega_k (n^2 \bs E_k, \bs G) + 
	\ell_0g_0\ga(n g_k \bs E_k, \bs G) &= 0,
	&\ \bs G \in \hHcurl,
\end{array}
\right.
\end{equation}
where \vspace*{-10pt}
\begin{align}
	\hat \trial_1 &:= \big\{\hat{\bs q} \in H^{-\frac{1}{2}}(\text{curl}, \Gamma_h): \bs n \times \hat{\bs q} = \bs n \times \bs E_{0,k} \text{ on } \Gamma \big\} , \\
	\hat \trial_2 &:= \big\{\hat{\bs q} \in H^{-\frac{1}{2}}(\text{curl}, \Gamma_h) \big\} ,
\end{align}
and $h$ denotes element-wise operations. The trace space on the mesh skeleton $\Gamma_h$ is defined as follows:
\begin{equation}
	H^{-\frac{1}{2}}(\text{curl}, \Gamma_h) :=
	\big\{
	\hat{\bs q} \in \hskip -2pt \prod_{K \in \Omega_h} H^{-\frac{1}{2}}(\text{curl}, \p K):
	\exists \bs q \in H(\text{curl}, \Omega):
	\gamma_t(\bs{q}|_K) = \hat{\bs q} 
	\big\} ,
	\label{eq:curl-trace-space}
\end{equation}
where the element boundary $\p K$ is assumed to be Lipschitz for every element $K \in \Omega_h$, and the continuous and surjective tangential trace operator is defined element-wise \cite{carstensen2016breaking},
\begin{equation}
	\gamma_t: H(\text{curl}, \Omega_h) \rightarrow
	\prod_{K \in \Omega_h} H^{-\frac{1}{2}}(\text{curl}, \p K) .
	\label{eq:tangent-trace}
\end{equation}
We employ the adjoint graph norm on the test space \cite{demkowicz2017dpg}.

\paragraph{Computational complexity}
The cross-sectional fiber geometry is approximated with four prismatic and twelve hexahedral isoparametric curvilinear elements, using orientation-embedded shape functions from \cite{fuentes2015shape}; in the direction of wave propagation, we use two elements per wavelength. For the numerical results presented here, we have not exploited mesh adaptivity; for this, we refer to \cite{henneking2020pollution} where the efficacy of the robust DPG error indicator for capturing higher-order modes in a waveguide is demonstrated. We use fast integration techniques to accelerate the element-wise computation of optimal test functions \cite{mora2019fast,badger2020fast}. In all numerical simulations, the (isotropic) order of approximation is $p=5$. The resulting DPG linear system has approximately 30,000 degrees of freedom per wavelength, i.e., ca.\ 3 million (complex double precision) unknowns for a fiber with 100 wavelengths. The weakly-coupled Maxwell systems are solved via Picard iteration which requires multiple solves per time step. However, the elongated (cylindric) geometry of the fiber enables efficient solves with a nested dissection elimination ordering. For the results presented here, we have used the MUMPS solver package \cite{amestoy2001mumps} for which a single solve of our model with 100 wavelengths requires ca.\ 2 minutes on a modern manycore compute architecture. All of our simulations were computed on Skylake (SKX) compute nodes using the Stampede2 supercomputer at the Texas Advanced Computing Center; OpenMP threading was enabled to exploit shared-memory parallelism. Preliminary tests have shown that the nested dissection ordering for the fiber amplifier model also scales well with MPI distributed memory computing on many compute nodes (some examples for linear waveguides with several thousand wavelengths are shown in \cite{henneking2020pollution}).

\begin{table}[htb]
	\caption{Step-index fiber: parameters}
	\centering
	\begin{tabular}{lll}
		Symbol & Description & Value \\
		\hline
		$r_{\text{core}}$ & Core radius & $12.7 \ \mu$m \\
		$r_{\text{clad}}$ & Cladding radius & $127 \ \mu$m \\
		$n_{\text{core}}$ & Refractive index in fiber core & $1.4512$ \\
		$n_{\text{clad}}$ & Refractive index in fiber cladding & $1.4500$ \\
		NA & Core numerical aperture & $0.059$ \\
		$\lambda_s$ & Signal wavelength & 1064 nm \\
		$\lambda_p$ & Pump wavelength & 976 nm \\
		$V(\omega_s)$ & Normalized signal frequency & $4.43$
	\end{tabular}
	\label{tab:fiber-parameters}
\end{table}

\begin{figure}[htb]
	\centering
	\begin{subfigure}[b]{0.48\textwidth}
		\includegraphics[width=\textwidth]{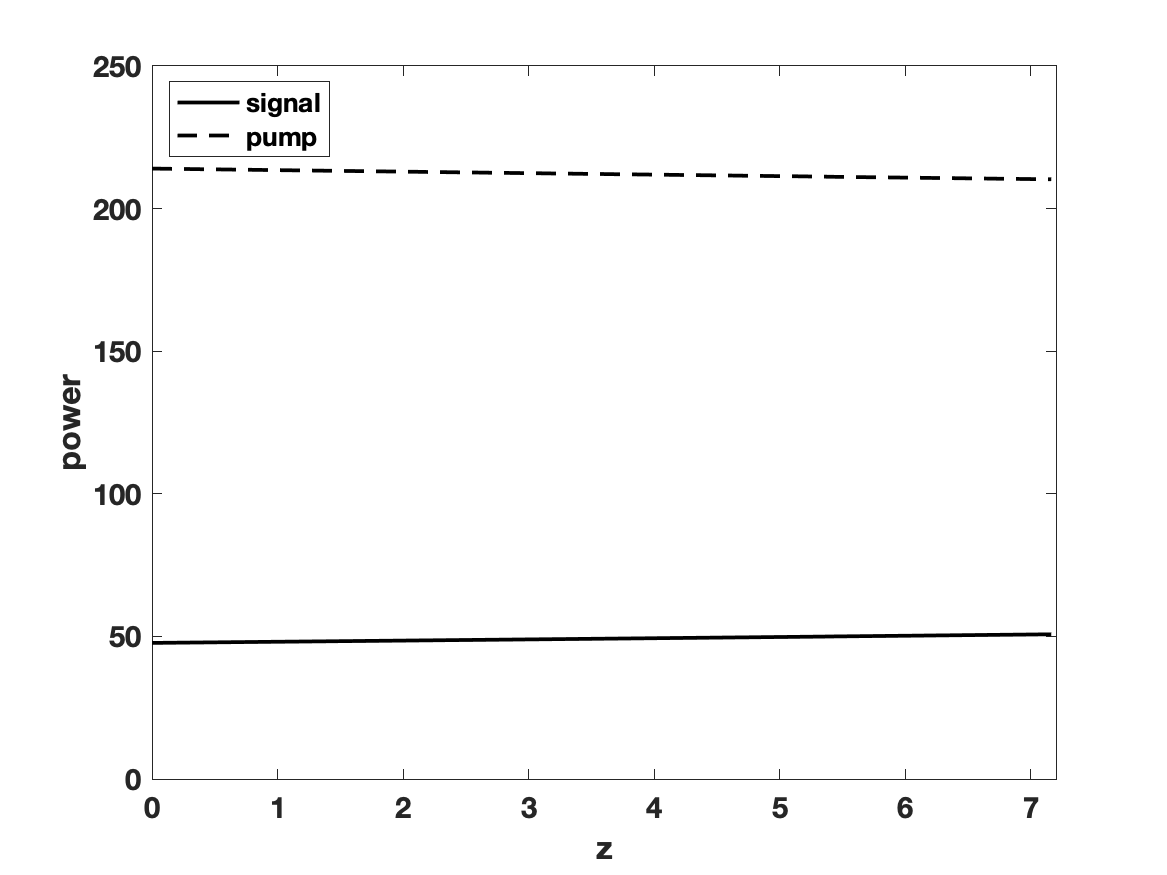}
		\caption{$\ga = 10^2$}
	\end{subfigure}
	\begin{subfigure}[b]{0.48\textwidth}
		\includegraphics[width=\textwidth]{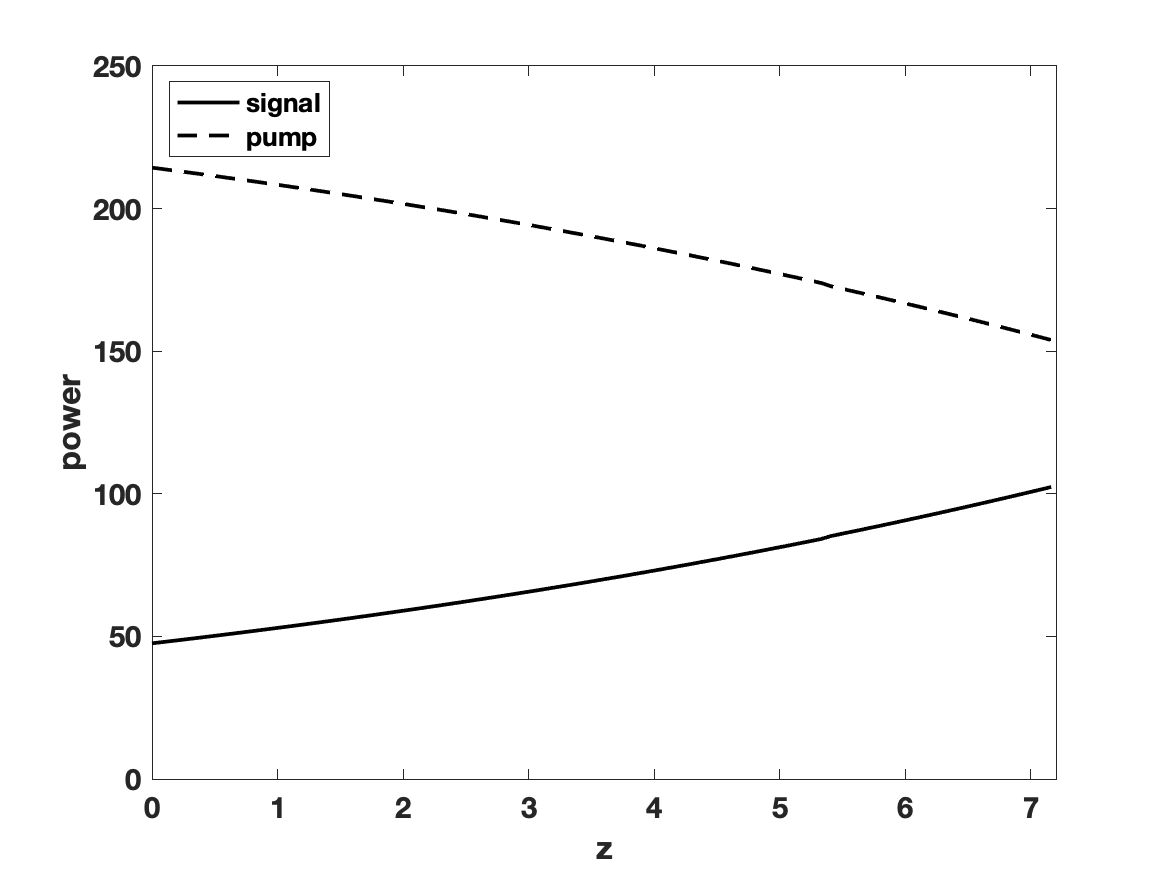}
		\caption{$\ga = 1.25 \cdot 10^3$}
	\end{subfigure}
	\begin{subfigure}[b]{0.48\textwidth}
		\includegraphics[width=\textwidth]{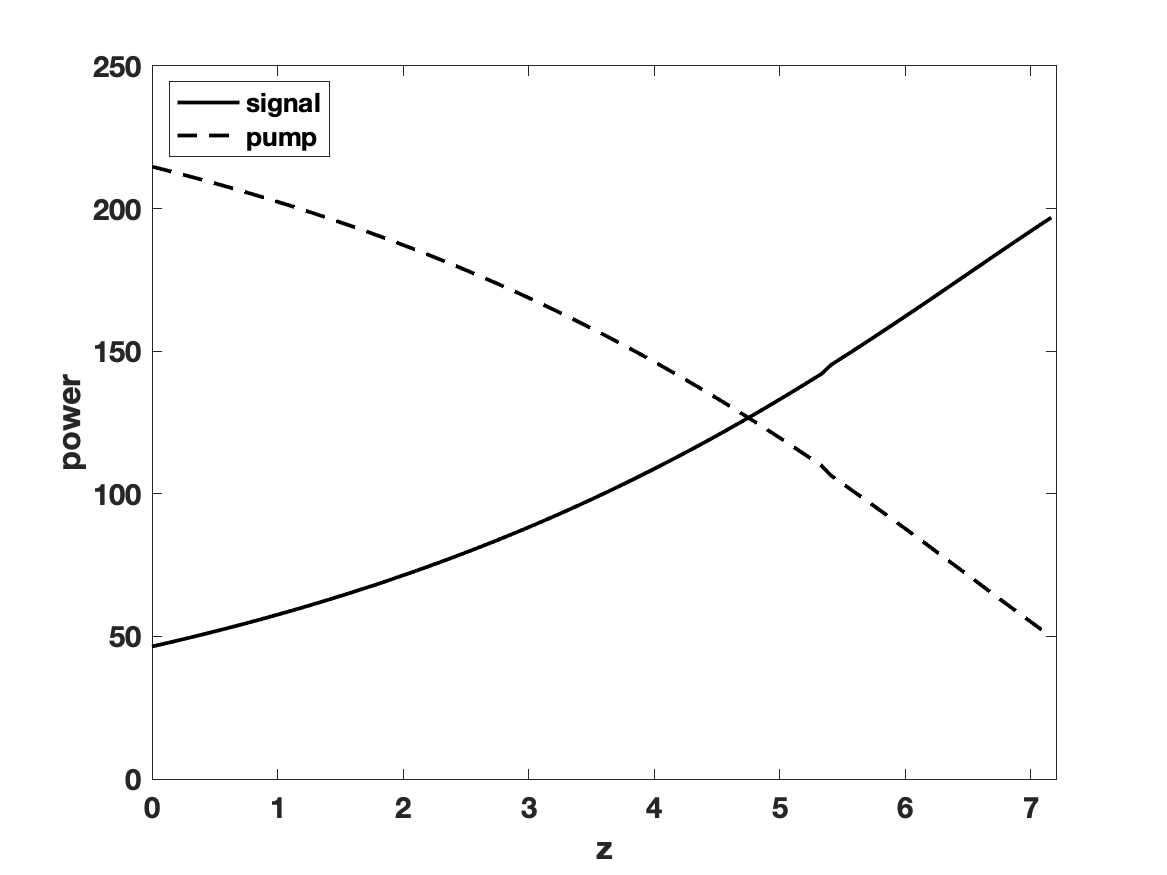}
		\caption{$\ga = 2.5 \cdot 10^3$}
	\end{subfigure}
	\begin{subfigure}[b]{0.48\textwidth}
		\includegraphics[width=\textwidth]{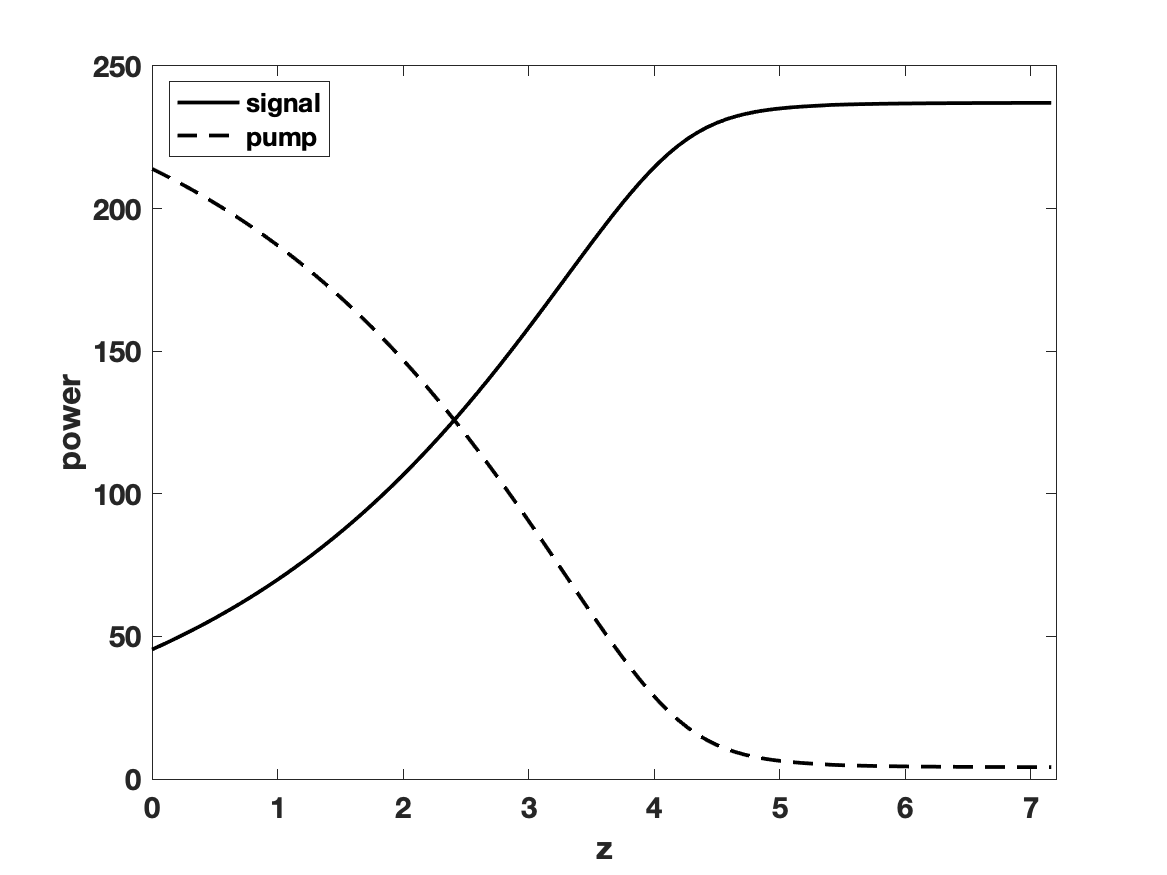}
		\caption{$\ga = 5 \cdot 10^3$}
	\end{subfigure}
	\caption{Active gain Maxwell simulation for 120 wavelengths: power distribution along the fiber amplifier in the signal and the pump field for different gain scaling terms $\ga$. Even within a short fiber, the entire pump energy can be transferred into the signal if the gain amplification is scaled artificially.}
	\label{fig:gain-amplifier}
\end{figure}

\paragraph{Nonlinear gain experiments}
The fiber parameters used in the numerical simulations of our model are given in Table~\ref{tab:fiber-parameters}. Results are shown for a core-pumped, co-pumped, fiber amplifier where the pump and signal fields are both excited with their respective fundamental mode at the fiber input. First, we investigate the cross-sectional power flux inside the fiber amplifier along the longitudinal axis. The power flux through a surface orthogonal to the fiber axis can be computed directly from the trace unknowns in the DPG broken ultraweak formulation \cite{nagaraj2018raman}. Figure~\ref{fig:gain-amplifier} illustrates the effect of the artificial gain amplification term $\ga$. Recall that the depicted non-dimensional values correspond to the appropriate dimensional scales (as shown in Table~\ref{tab:scales}). For a fiber of about 120 wavelengths, a small (relative to the ratio $L/\tilde L$) value of $\ga=10^2$ causes almost no exchange of power between pump and signal fields. For $\ga = 5 \cdot 10^3$, the gain term suffices to transfer the power of the pump into the signal within the short fiber.

\begin{algorithm}[htb]
\begin{algorithmic}
	\STATE $\bs E_{s}^{(0)} \gets \bs 0$
	\STATE $\bs E_{p}^{(0)} \gets \bs 0$
	\FOR{$i=1$ \TO $i_{max}$}
		\STATE $\bs E_{s}^{(i)} \gets$ solve with $g_{s} := g_{s} \big(\bs E_{s}^{(i-1)}, \bs E_{p}^{(i-1)} \big)$
		\STATE $\bs E_{p}^{(i)} \gets$ solve with $g_{p} := g_{p} \big(\bs E_{s}^{(i-1)}, \bs E_{p}^{(i-1)} \big)$
		\IF {$\| \bs E_{s}^{(i)} - \bs E_{s}^{(i-1)} \| / \| \bs E_{s}^{(i)} \| < \varepsilon$}
			\STATE \textbf{break}
		\ENDIF
	\ENDFOR
\end{algorithmic}
	\caption{Nonlinear gain problem: Picard iteration}
	\label{alg:Picard}
\end{algorithm}

\begin{figure}[htb]
	\centering
	\begin{subfigure}[b]{0.48\textwidth}
		\includegraphics[width=\textwidth]{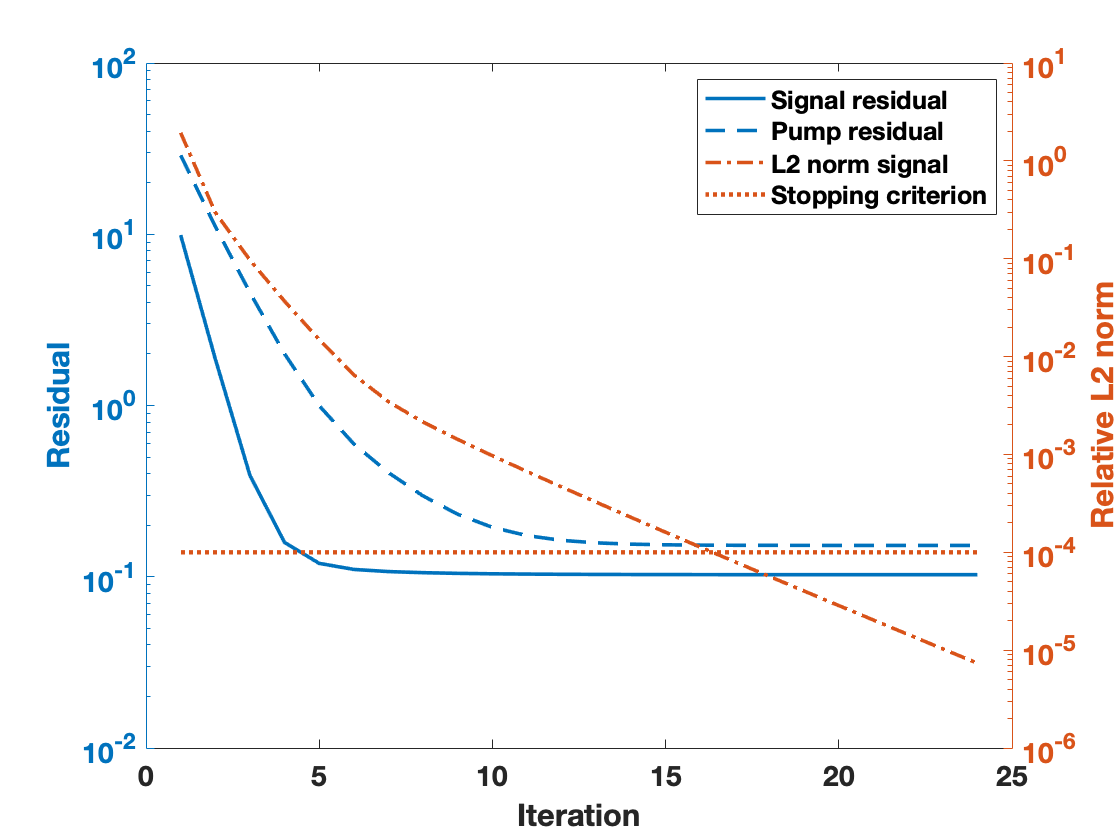}
		\caption{Residual in nonlinear solve}
		\label{fig:gain-residual}
	\end{subfigure}
	\begin{subfigure}[b]{0.48\textwidth}
		\includegraphics[width=\textwidth]{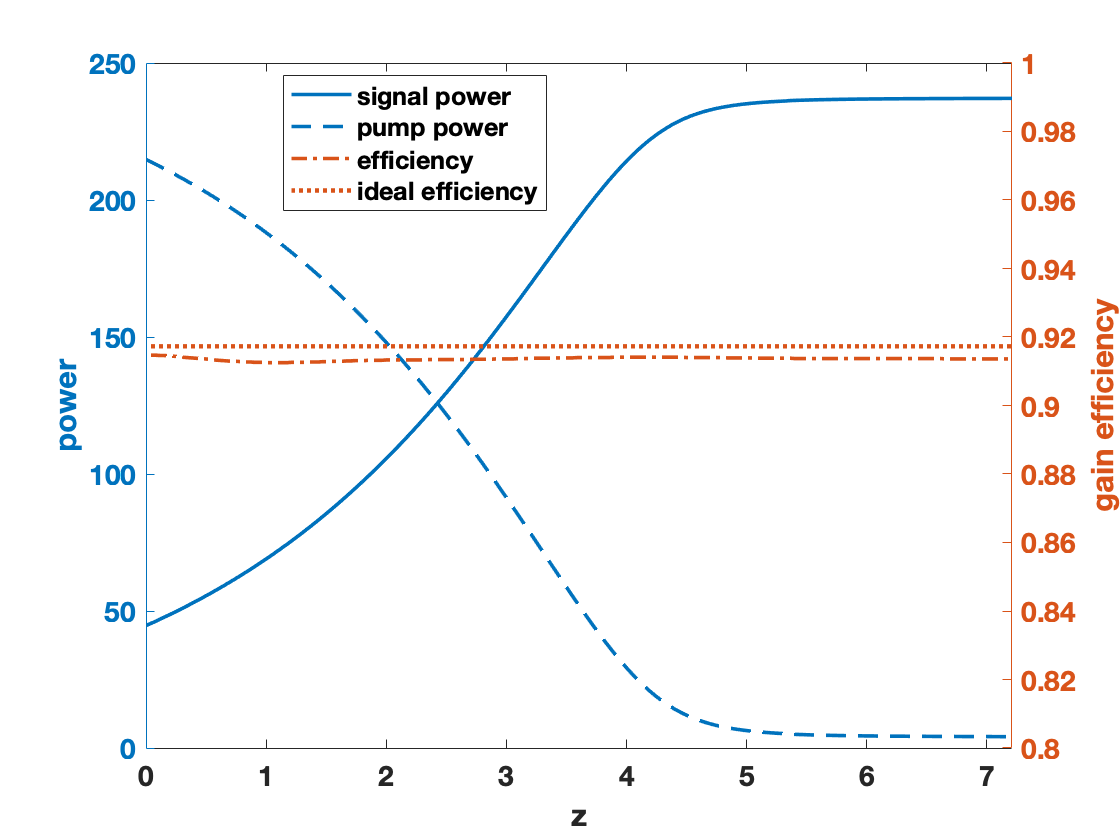}
		\caption{Amplifier efficiency}
		\label{fig:gain-efficiency}
	\end{subfigure}
	\caption{Active gain Maxwell simulation for $120$ wavelengths: (a) the DPG residual, measured in the adjoint graph norm, converges in the Picard iterations; and (b) the efficiency of the power transfer from the pump into the signal field is near the ideal amplifier efficiency.}
\end{figure}

The stopping criterion for the nonlinear solve is based on the relative change of the signal field measured in the $L^2$ norm. The implementation of the Picard iteration is given in Algorithm~\ref{alg:Picard}. Alternatively, one could use the DPG residual directly as a stopping criterion, but this would add significant computational cost in every iteration. Preferably, the residual should only be computed at the very end of each nonlinear solve in order to determine which elements to refine in an adaptive mesh refinement. The convergence of the residual for each field is plotted alongside the convergence of the relative signal field, measured in $L^2$, and its stopping criterion ($\eps = 10^{-4}$) in Figure~\ref{fig:gain-residual}. One validation of the model itself is presented in Figure~\ref{fig:gain-efficiency}, showing the obtained pump efficiency at different points along the fiber, compared to the ideal efficiency $\lambda_p/\lambda_s \approx 91.7\%$. As expected, the optical-to-optical efficiency in the fiber is slightly below this ideal efficiency since the light guided in the core falls into transverse modes which do not uniformly overlap the core region, which means that the gain is not perfectly/maximally saturated.

\begin{figure}[ht]
	\centering
	\begin{subfigure}[b]{0.48\textwidth}
		\includegraphics[width=\textwidth]{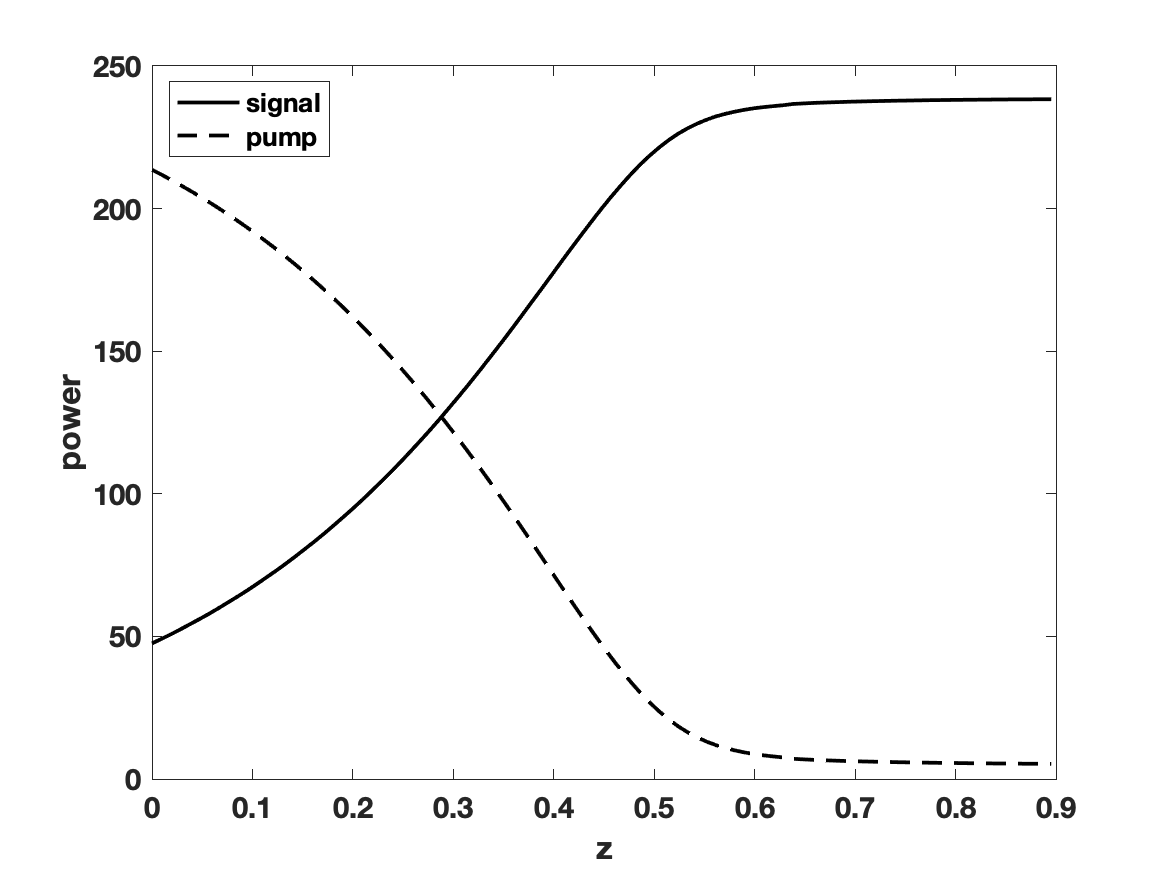}
		\caption{$15\lambda,\ \ga = 4 \cdot 10^4$}
		\label{fig:gain-scaling-15}
	\end{subfigure}
	\begin{subfigure}[b]{0.48\textwidth}
		\includegraphics[width=\textwidth]{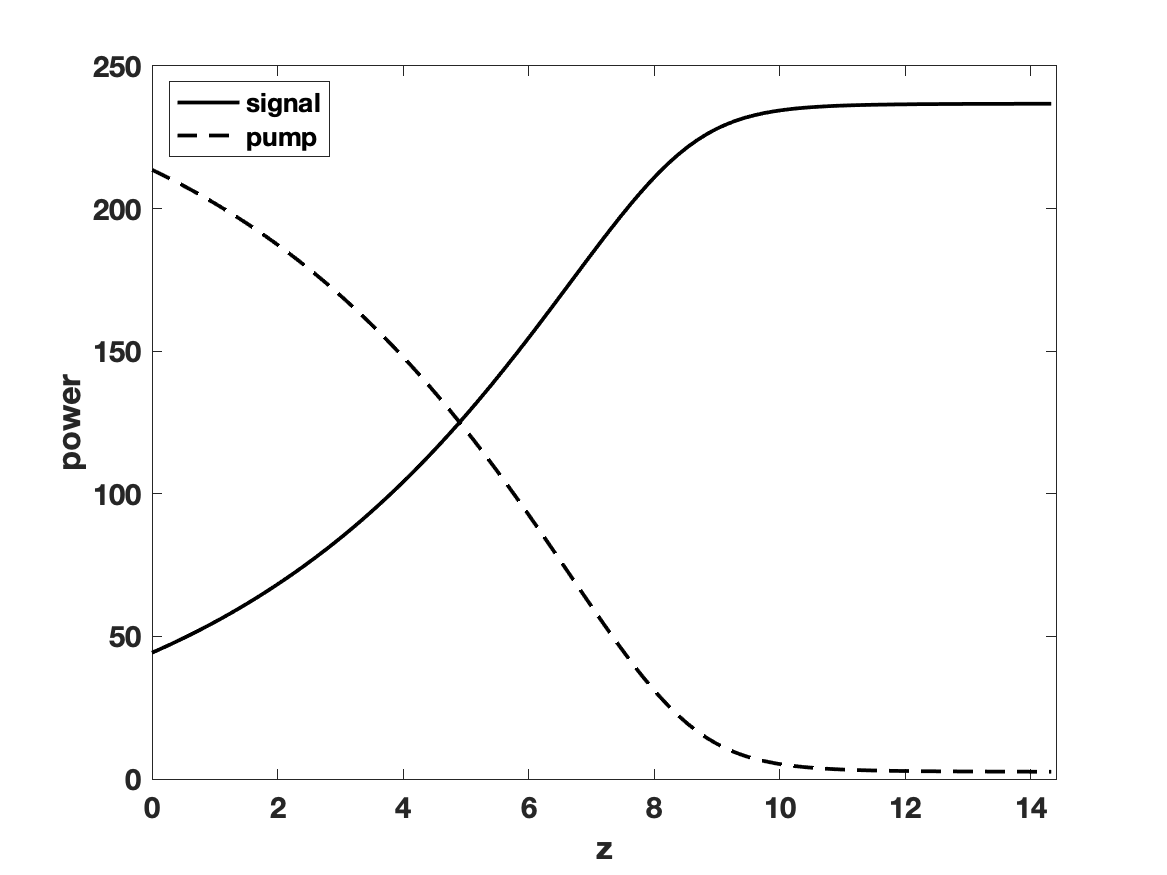}
		\caption{$240\lambda,\ \ga = 2.5 \cdot 10^3$}
		\label{fig:gain-scaling-120}
	\end{subfigure}
	\caption{Active gain Maxwell simulation: scaling experiment. Approximately the same power distribution (for both signal and pump field) is obtained for two fibers of different length by scaling the gain term $\ga$ appropriately.}
	\label{fig:gain-scaling}
\end{figure}

\paragraph{Gain scaling experiment}
When we introduced the artificial gain scaling term $\ga$, it was with the goal to preserve certain quantities of interest in the scaled model. A numerical verification for the preservation of the power distribution in fibers of different lengths is shown in Figure~\ref{fig:gain-scaling}. We compare a fiber of 240 wavelengths (the largest one computed in our experiments) with $\ga = 2.5 \cdot 10^3$ to a shorter fiber of 15 wavelengths and $\ga = 4 \cdot 10^4$ (appropriately scaled by a factor of 16).

Table~\ref{tab:iterations} summarizes the number of iterations needed in the nonlinear solve until the stopping criterion is reached. As $\ga$ is increased for a fiber of fixed length, the number of iterations increases, indicating the stronger nonlinearity of the system---previously illustrated by Figure~\ref{fig:gain-amplifier}. On the diagonals of the table (starting from top right towards the bottom left) the number of iterations remains constant. These diagonals represent the ``gain scaling'' of a fiber: moving along them means increasing the fiber length and decreasing $\ga$ by the same factor.

\begin{table}[ht]
	\caption{Active gain Maxwell simulation: number of iterations until convergence criterion is satisfied for different number of wavelengths ($\# \lambda$) and gain scaling terms ($\ga$). More Picard iterations are needed for stronger gain amplification and/or longer fibers.}
	\label{tab:iterations}
	\centering
	\begin{tabular}{c|c:c:c:c:c:c}
			\diagbox[height=1.5\line]{$\# \lambda$}{$\ga$} & $1.25 \cdot 10^3$ & $2.5 \cdot 10^3$ & $5 \cdot 10^3$ & $10^4$ & $2 \cdot 10^4$ & $4 \cdot 10^4$ \\
			\hline
			15 & \tc{black}{3} & \tc{black}{4} & \tc{black}{4} & \tc{black}{4} & \tc{black}{8} & \tc{black}{15}\\
			30 & \tc{black}{4} & \tc{black}{4} & \tc{black}{5} & \tc{black}{9} & \tc{black}{16} & -\\
			60 & \tc{black}{4} & \tc{black}{5} & \tc{black}{10} & \tc{black}{16} & - & -\\
			120 & \tc{black}{5} & \tc{black}{10} & \tc{black}{15} & - & - & -
	\end{tabular}
\end{table}

While these numerical verifications do not rigorously prove that scaled short fibers reproduce accurate power distributions for real-length fibers where the scaling factor is much larger, they do indicate that the fiber may at least be scaled to some extent in our model while preserving the signal and pump power curves with sufficient accuracy.

\paragraph{DPG primal heat}
After solving the initial Maxwell problem at ambient temperature, time is advanced via the time-stepping of the heat equation. The operator form of the heat problem is,
\begin{equation}
\left\{
\begin{array}{rll}
	\frac{\p (\dT)}{\p t} - 
	\alpha_0 \nabla \cdot ( \Lz \nabla (\dT) )
	&= Q_0 Q(I_{\{s,p\}}) & \text{in } \Omega \times (0,t_{\text{max}}], \\
	\dT &= 0 & \text{on } \Gamma \times [0,t_{\text{max}}], \\
	\dT &= 0  & \text{in } \Omega \times \{ 0 \}.
\end{array}
\right.
\label{eq:heat-operator}
\end{equation}

The DPG broken primal formulation with implicit Euler time-stepping is,
\[
\left\{
\begin{array}{l}
	\dT_{n+1} \in H_0^1(\Omega),\
	\hat{\bs \sigma}_{n+1} \cdot \bs n \in H^{-\frac{1}{2}}(\Gamma_h), \\[3pt]
	(\dT_{n+1}, v) + 
	\dt \hskip 1pt \alpha_0 (\Lz \nabla (\dT_{n+1}), \nabla v)
	- \dt \hskip 1pt \alpha_0 \lb \hat{\bs \sigma}_{n+1} \cdot \bs n, v \rr_{\Gamma_h}
	= (\dT_n,v) + \dt \hskip 1pt Q_0 (Q(I_{\{s,p\}}),v),\
	v \in H^1(\Omega_h),
\end{array}
\right.
\]
where $n = 0,1,\ldots, N-1,$ and $\delta t = t_{\text{max}} / N$. The primal formulation includes the heat flux as an unknown on the mesh skeleton, where the trace space $H^{-\frac{1}{2}}(\Gamma_h)$ can be introduced as a product space of element-wise normal traces, similar to (\ref{eq:curl-trace-space})--(\ref{eq:tangent-trace}). We equip the test space with the standard energy norm. With this setup, the coupled Maxwell/heat problem now requires the use of elements of the entire $H^1-H(\text{curl})-H(\text{div})-L^2$ exact sequence \cite{carstensen2016breaking,fuentes2015shape}. 

\begin{figure}[htb]
	\centering
	\includegraphics[width=0.8\textwidth]{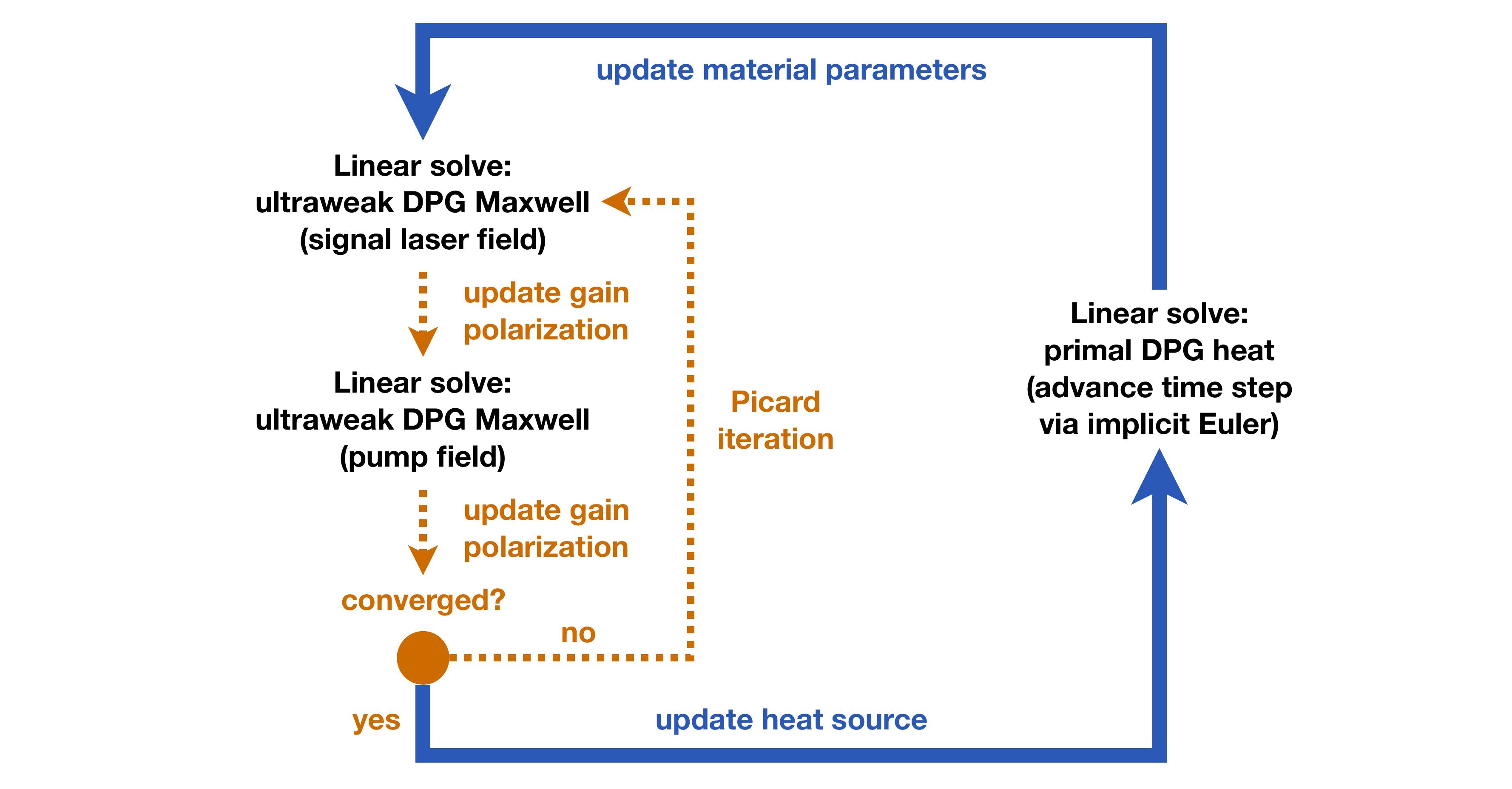}
	\caption{Coupled Maxwell/heat simulation: overview of the numerical scheme. Picard iterations are used for linearizing the weakly-coupled Maxwell systems, and implicit Euler time-stepping for the heat equation advances the time step of the simulation. Material parameters and heat deposition in the fiber are updated once per time step.}
	\label{fig:coupled-scheme}
\end{figure}

Figure \ref{fig:coupled-scheme} displays the overall numerical scheme for the coupled problem. The weakly-coupled Maxwell systems are solved via Picard iteration, where the gain polarization is updated once per iteration (as shown in Algorithm~\ref{alg:Picard}), and the simulation's time step is advanced via implicit Euler time-stepping in the heat equation. The fiber material parameters (thermally induced refractive index perturbation) and heat source (heat deposition in the fiber) are updated once per time step.

\vspace{-5pt}

\begin{figure}[htb]
	\centering
	\begin{subfigure}[b]{0.48\textwidth}
		\includegraphics[width=\textwidth]{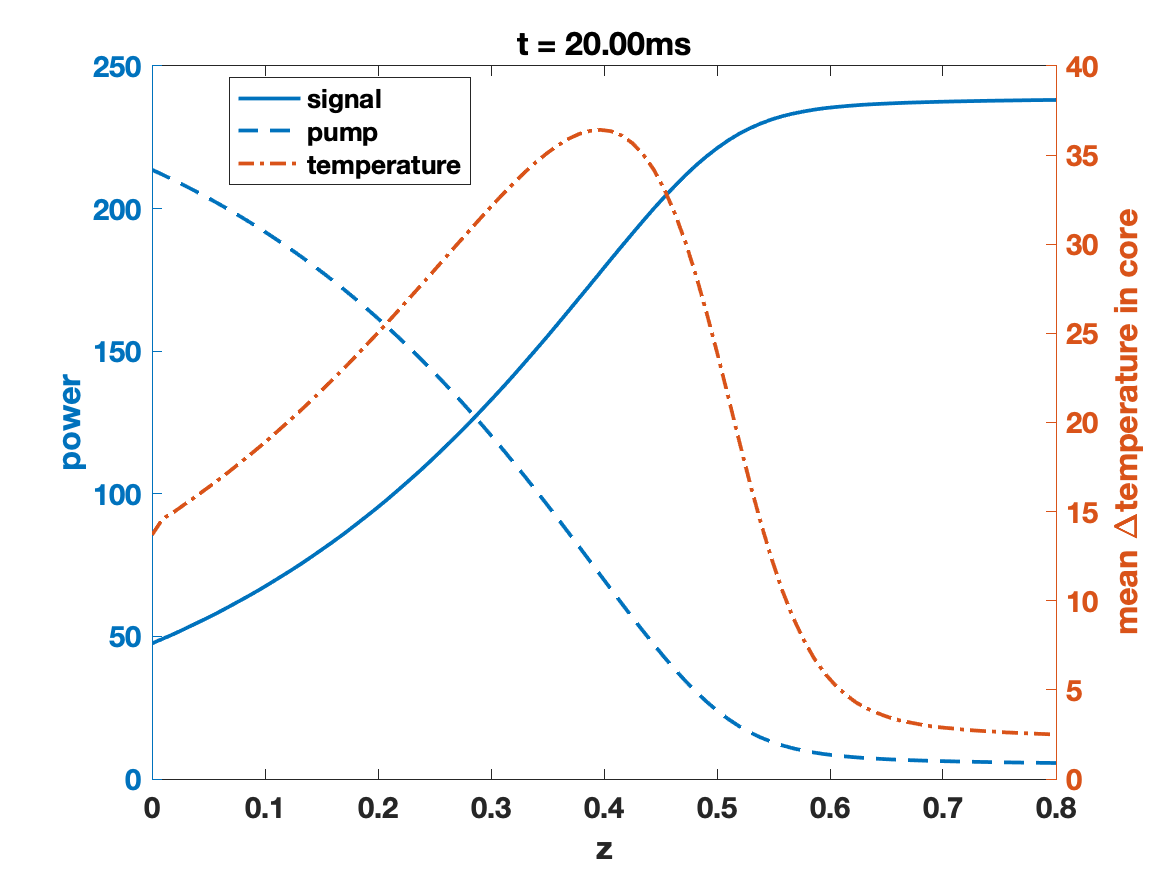}
		\caption{$15\lambda,\ \ga = 4 \cdot 10^4,\ \az = 1.2 \cdot 10^{-6}$}
		\label{fig:temp-scaling-15}
	\end{subfigure}
	\begin{subfigure}[b]{0.48\textwidth}
		\includegraphics[width=\textwidth]{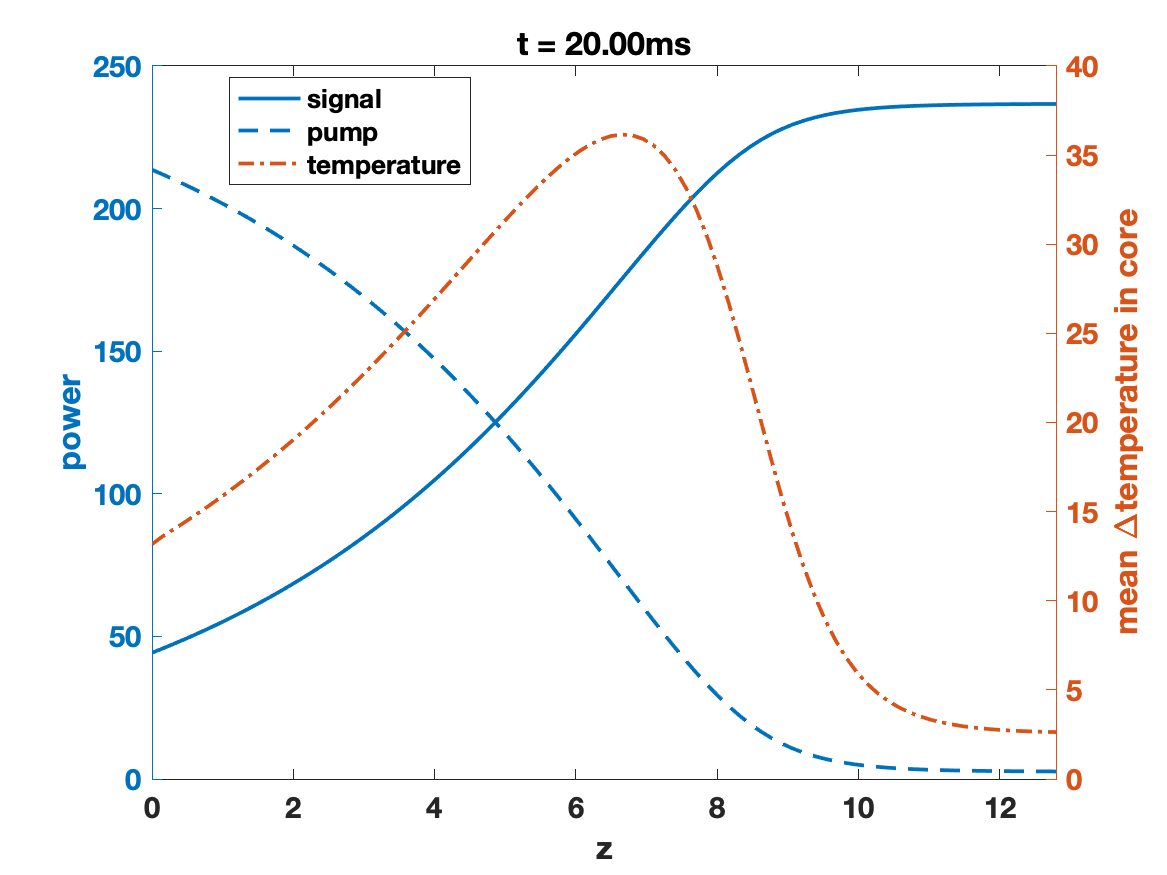}
		\caption{$240\lambda,\ \ga = 2.5 \cdot 10^3,\ \az = 1.92 \cdot 10^{-5}$}
		\label{fig:temp-scaling-240}
	\end{subfigure}
	\caption{Coupled Maxwell/heat simulation: scaling experiment. Approximately the same temperature distribution is obtained after $200$ time steps ($\delta t = 0.1$ ms) for two fibers of different length by scaling the anisotropic heat diffusion coefficient $\az$ appropriately.}
	\label{fig:temp-scaling}
\end{figure}

\paragraph{Heat scaling experiment and induced thermal perturbations}
In the numerical experiments, $\dt=0.1$ ms has proven to be sufficiently small to provide stability in the heat solve. For a given solution to the nonlinear Maxwell problem, the transient heat equation attains steady-state after circa $15$ ms. Figure~\ref{fig:temp-scaling} shows the temperature distribution along the fiber amplifier after $20$ ms for two different fibers: one with 15 wavelengths (Figure~\ref{fig:temp-scaling-15}) and one with 240 wavelengths (Figure~\ref{fig:temp-scaling-240}). The results indicate that the heat distribution in the fiber can be computed accurately on a short fiber, as expected by the anisotropic diffusion operator scaling argument. In other words, to obtain an accurate heat curve in the scaled model, one only requires the accurate power distribution (or intensity) of the signal and pump field along the fiber.

\begin{figure}[htb]
	\centering
	\includegraphics[width=0.96\textwidth]{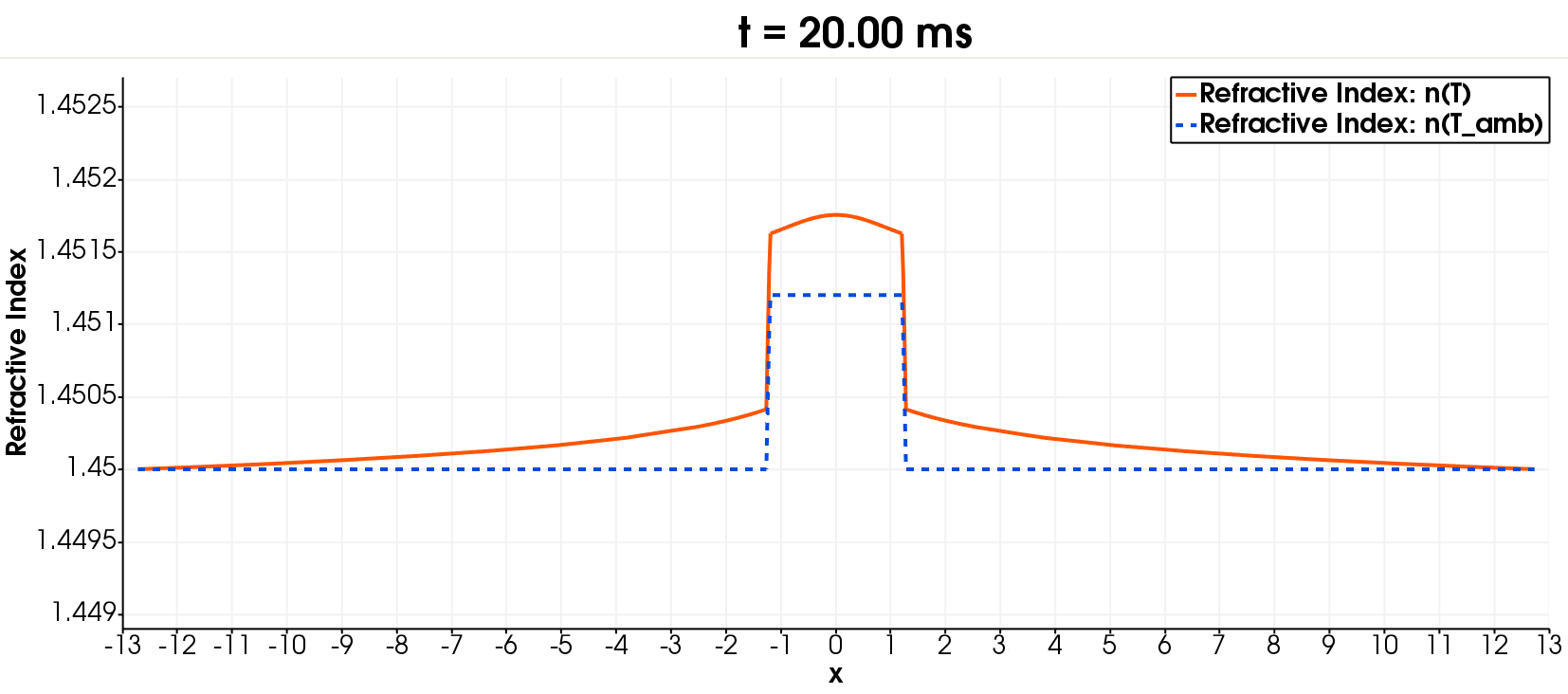}
	\caption{Coupled Maxwell/heat simulation: refractive index profile inside the fiber amplifier, plotted along the $x$-axis in a slice orthogonal to the longitudinal fiber axis. Near the regions with peak temperature, we observe a significant thermally induced perturbation compared to the step-index profile at ambient temperature.}
	\label{fig:ref-index}
\end{figure}

To investigate the nonlinear effects induced by thermal perturbations in the fiber amplifier, we compute the solution to the weakly-coupled Maxwell systems in every time step of the heat equation with updated material parameters. The thermally induced perturbations are represented by the change of the material refractive index inside the fiber core and cladding. Figure~\ref{fig:ref-index} shows the refractive index plotted along $x$ across a slice of the fiber orthogonal to the longitudinal axis close to where the peak temperature occurs. It shows that the refractive index profile is perturbed significantly compared to the step-index profile at ambient temperature. Consequently, the guided propagating fields may be perturbed in a significant way as the temperature develops inside the fiber amplifier, e.g., the core-guided fields will likely experience thermal lensing.

\section{Summary}
\label{sec:summary}
%
%

Our nonlinear fiber amplifier model has now been augmented with active gain and an integrated thermal response.
The simulation incorporates two weakly-coupled time-harmonic 3D vectorial Maxwell systems for the propagating electromagnetic signal and pump fields.
This coupling between the fields occurs both through the active gain in the fiber core region and through the thermally induced refractive index perturbations that result from the heat deposition caused by the lasing (see Figure~\ref{fig:ref-index}).
The steady-state ytterbium ion population concentrations are updated according to how the irradiances of the pump and signal fields evolve along the length of the fiber. As expected, the pump field experiences loss (negative gain), and the laser signal experiences amplification (positive gain) such that the ideal optical-to-optical efficiency is never surpassed (see Figure~\ref{fig:gain-efficiency}). Furthermore, this gain via stimulated emission is ultimately what drives the heating in the fiber, resulting in a peak heat load in the region of greatest transfer of energy from the pump field to the signal field (see Figure~\ref{fig:temp-scaling}).

The proposed model differs from the more common BPM or CMT approaches for numerical simulation of fiber lasers in many ways. Most importantly, our model makes very few assumptions on the propagating fields in order to provide a high-fidelity simulation tool. This approach leads to high computational cost when solving the fiber for many wavelengths in a high-order finite element discretization. However, the resulting linear system can be efficiently solved with a nested dissection elimination ordering. The DPG method is used for discretizing both the Maxwell systems and the heat equation, yielding a stable discretization with a built-in error indicator. The data from the numerical solution of the coupled Maxwell/heat 3D fiber amplifier model enables the analysis of the interplay between thermal perturbations of the material and the propagating electromagnetic fields with great accuracy. Additionally, the generality of the model makes it possible to add further nonlinearities, compute counter- or bi-directional pumping configurations, or study birefringent (anisotropic) fibers, all of which are difficult to realize in simpler fiber models due to their inherent assumptions. At the same time, the high accuracy of our model can be exploited to validate typical approximations made in lower-fidelity models.

In order to make the computation feasible for a full-length fiber, we introduced a longitudinal scaling for the Maxwell and heat problem. This was accomplished with an artificial material parameter that effectively enhances the gain per unit length, which was shown to preserve certain quantities of interest, e.g., amplifier efficiency, within the tested parameter regime. Of particular importance, we were able to scale the heat equation naturally through an anisotropic diffusion operator obtained by a change of coordinates such that the only source of artificial error introduced into the amplifier model by this scaling is through the gain along the fiber.

In summary, we were able to make significant modeling advancements, demonstrating how to incorporate nonlinear laser amplification and thermal effects into the full vectorial Maxwell model for an ytterbium-doped fiber amplifier, while maintaining the computational feasibility. We believe that this unique model effectively complements other, typically lower-fidelity, fiber amplifier models.

In future work, we plan to leverage this model to investigate mode instabilities; it is expected that the artificially scaled fiber needs to be long enough to capture the interference patterns between the transverse guided core modes, meaning at least multiple mode beat lengths. For that reason, a distributed finite element code for this model is currently under development. Our preliminary results for linear waveguides indicate that the simulation of our fiber amplifier model with at least 10,000 wavelengths is feasible with current compute architectures. More computational power would also make it feasible to model more complex fiber geometries, e.g., microstructure fibers, or study the effect of coiling the fiber around a cooling spool.

\FloatBarrier
\section*{Acknowledgement}
This work was supported by AFOSR grant no.\ FA9550-19-1-0237.

\bibliographystyle{abbrv}

\bibliography{./shortref,./ref,./ref_laser}

\end{document}